\input amstex
\magnification=\magstep1
\documentstyle{amsppt}
\nologo
\NoBlackBoxes
\pageheight{8.5truein}
\pagewidth{6.5truein}
\def\var{\varepsilon}

\define\fc{\frac 1 2}
\define\sumstar{\operatornamewithlimits{\sum\nolimits^\star}}
\define\ssum{\operatornamewithlimits{\sum\,\sum}}
\define\sn{\smallskip\noindent}
\define\mn{\medskip\noindent}
\define\bn{\bigskip\noindent}

\define\fk{\frak a}

\define \nl{\newline}
\documentstyle {amsppt}
\topmatter
\title\nofrills{Spacing of Zeros of Hecke $L$-Functions \\
and the Class Number Problem} \endtitle
\rightheadtext\nofrills{Spacing of Zeros of Hecke $L$-Functions and
the Class Number Problem}
\leftheadtext\nofrills{B. Conrey and H. Iwaniec}
\author {B. Conrey$^*$ and H. Iwaniec$^{**}$
\thanks{\null\newline
${}^*$Supported by NSF grant DMS - 99-80392 \null\newline
${}^{**}$Supported by NSF grant DMS-98-01642  and FRG grant DMS-00-74028} \endthanks} \endauthor
\address   Brian Conrey \endgraf
American Institute of Mathematics \endgraf
360 Portage Avenue \endgraf
Palo Alto, CA  94306-2244 \endgraf
USA \endgraf
{\it E-mail address:} {\bf conrey\@best.com}
\endaddress

\address
Henryk Iwaniec \endgraf
Mathematics Department \endgraf
Rutgers University \endgraf
Hill Center \endgraf
Piscataway, NJ  08854-8019 \endgraf
USA \endgraf
{\it E-mail address:} {\bf iwaniec\@math.rutgers.edu}
\endaddress

\endtopmatter
\document

{\bf Contents:}
\sn
1. A bit of history and results \nl
2. Basic automorphic forms \nl
3. Summation formulas \nl
4. Convolution sums \nl
5. Point to integral mean-values of Dirichlet's series \nl
6. Evaluation of ${\Cal A}(T)$ \nl
7. Approximate functional equation \nl
8. Evaluation of $\ell(s)$ on average \nl
9. Estimation of $x(s)$ on average \nl
10. Applications
\newpage

\subhead  1.  A Bit of History and Results\endsubhead
\medskip
The group of ideal classes ${\Cal C}\ell(K)$ of an imaginary
quadratic field $K = \Bbb Q(\sqrt {-q})$ is the most fascinating
finite group in arithmetic.  Here we are faced with one of the
most challenging problems in analytic number theory,  that is to
estimate the order of the group $h = |{\Cal C}\ell(K)|$.  C.F.
Gauss conjectured (in a parallel setting of binary quadratic
forms) that the class number $h = h(-q)$ tends to infinity as $-q$
runs over the negative discriminants.  Hence there are only a
finite number of imaginary quadratic fields with a given class
number. But how many of these fields are there exactly for $h =
1$, or $h = 2$, etc. ?  To answer this question one needs on
effective lower bound for $h$ in terms of $q$ (a fast computer
could be helpful as well).
\medskip
The problem was linked early on to the L-series $$ L(s,\chi) =
\sum^\infty_{n=1} \chi(n) \,n^{-s}  \tag 1.1 $$ for the real
character $\chi$ of conductor $q$ (the Kronecker symbol) $$
\chi(n) = \biggl( \frac{-q}{n} \biggr).  \tag 1.2 $$ In this
connection L. Dirichlet established the formula $$ h =
\pi^{-1}\sqrt q L(1, \chi)  \tag 1.3 $$ (we assume that $-q$ is a
fundamental discriminant, $q>4$, so there are two units $\pm 1$ in
the ring of integers ${\Cal O}_K \subset K$).  Rather than
estimating the class number, Dirichlet inferred from (1.3) that
$L(1,\chi)$ does not vanish, which property he needed to establish
the equidistribution of primes in arithmetic progressions.  Truly
the lower bound $$ L(1,\chi) \geqslant \frac{\pi}{\sqrt q}  \tag
1.4 $$ follows from (1.3), because $h \geqslant 1$.
\medskip
The Grand Riemann Hypothesis for $L(s,\chi)$ implies $$ (\log\log
q)^{-1} \ll L(1,\chi) \ll \,\log\log q,  \tag 1.5 $$ whence the
class number varies only slightly about $\sqrt q $ $$ \frac{\sqrt
q}{\log\log q} \ll h \ll \sqrt q \,\log\log q. \tag1.6 $$ But
sadly enough we may not see proofs of such estimates (which are
best possible in order of magnitude) in the near future.
\medskip
At present we know (after J. Hadamard and C.J. de la
Vall\'ee-Poussin) that $L(s,\chi) \neq 0$ for $s = \sigma + it$ in
the region $$ \sigma > 1 - \frac{c}{\log q (|t|+1)}  \tag 1.7 $$
where $c$ is a positive constant, for any character $\chi$(mod
$q$) with at most one exception.  The exceptional character
$\chi$(mod $q$) is real and the exceptional zero of $L(s,\chi)$ in
the region (1.7) is real and simple, say $\beta$ if it exists $$
\beta > 1- \frac{c}{\log q}.  \tag 1.8 $$ Using complex function
theory one can translate various zero-free regions of $L(s,\chi)$
which are stretched along the line Re $s =1$ to lower bounds for
$|L(s,\chi)|$ on the line Re $s =1$.  In the case of a real
character (1.2) H. Hecke (see [L1]) showed that if $L(s,\chi)$ has
no exceptional zero, then $L(1,\chi) \gg (\log q)^{-1}$, whence $$
h \gg \sqrt q \,(\log q)^{-1}  \tag 1.9 $$ by the Dirichlet
formula (1.3). Moreover if $L(s,\chi)$ does have an exceptional
zero $s = \beta$ satisfying (1.8), then we have quite precise
relations between $\beta$ and $h$ (see [GSc], [G1], [GS]).  In
particular one can derive from the Dirichlet estimate (1.4)that $$
\beta \leqslant 1- \frac{c}{\sqrt q}.  \tag 1.10 $$
\medskip
Back to the history we should point out that E. Landau [L1] first
come up with ideas which pushed the exceptional zero further to
the left of (1.10). Generalizing slightly in this context we owe
to Landau the product (a quadratic lift $L$-function) $$ Lan(s,f)
= L(s,f)L(s,f \otimes \chi) = \sum_n a_f(n) \,n^{-s} \tag 1.11 $$
where $$ L(s,f) = \sum_n \lambda_f(n) \,n^{-s} $$ can be any
decent L-function and $L(s,f \otimes\chi)$ is derived from
$L(s,f)$ by twisting its coefficients $\lambda_f(n)$ with
$\chi(n)$.  If $L(s,f)$ does have an Euler product so do $L(s,f
\otimes\chi)$ and $Lan(s,f)$. The key point is that the prime
coefficients of $Lan(s,f)$ are $$ a_f(p) = \lambda_f(p)
(1+\chi(p)). $$ Assuming the class number $h$ is small (or
equivalently that $L(s,\chi)$ has an exceptional zero) we find
that $\chi(p)=-1$ and $a_f(p) = 0$ quite often if $p \ll \sqrt q$.
In other words $\chi(m)$ pretends to be the M\"obius function
$\mu(m)$ on squarefree numbers.  Therefore, under this ficticious
assumption, $L(s,f\otimes\chi)$ approximates to $L(s,f)^{-1}$ and
$Lan(s,f)$ behaves like a constant (no matter what  $s$ and $f$
are \,\, !).
\medskip
Landau worked with $Lan(s,\chi')=
L(s,\chi^\prime)\,L(s,\chi\chi^\prime)$ where $\chi'$(mod $q'$) is
any real primitive character other than $\chi$(mod $q$). He [L1]
proved that for any real zeros $\beta, \beta^\prime$ of
$L(s,\chi), \,L(s,\chi^\prime)$ respectively $$ \min\,(\beta,
\beta^\prime) \leqslant 1- \frac{c}{\log q \,q^\prime}.  \tag 1.12
$$ This shows that the exceptional zero occurs very rarely (if at
all\,\,?).
\medskip
Next a repulsion property of the exceptional zero was discovered,
notably in the works by M. Deuring [D] and H. Heilbronn [H].  This
says - the closer $\beta$ is to the point $s=1$ the further away
from $s=1$ are the other zeros - not only of $L(s,\chi)$, but of
any L-function for a character of comparable conductor.  The power
of repulsion is masterly exploited in the celebrated work of Yu.V.
Linnik [L] on the least prime in an arithmetic progression.
\medskip
A cute logical play with repulsion led E. Landau [L2] to the lower
bound $$ h \gg q^{\frac 1 8 -\var}  \tag 1.13 $$ for any $\var >
0$, the implied constant depending on $\var$. Slightly later by
the same logic, but with more precise estimates for relevant
series  C.L. Siegel [S] proved that $$ h \gg q^{\frac 1 2 -\var}.
\tag 1.14 $$ Both estimates suffer from the serious defect of
having the implied constant not computable.  For that reason the
Landau-Siegel estimates do not help to determine all quadratic
imaginary fields with a fixed class number.  The case $h=1$ was
eventually solved by arithmetical means (complex multiplication)
by H. Heegner [He] and H. Stark[S] and by transcendental means
(linear forms in logarithms) by A. Baker [B] (see also the notes
[S2] about earlier attempts by A.O. Gelfond and Yu.V. Linnik
[GL]).
\medskip
By way of the repulsion one may still hope to produce effective
results provided an ``exceptional'' zero is given numerically.
But, believing in GRH one cannot expect to find a real zero of any
decent L-function other than at the central point $s=\frac 1 2$.
Hence the question: Does the central zero have an effect on the
class number? Yes it does, and the impact depends on the order of
the zero. This effect was first revealed in conversations by J.
Friedlander in the early 70's. Soon after J.V. Armitage gave an
example of the zeta function of a number field which vanishes at
the central point, Friedlander [F] succeeded in estimating
effectively the class number of relative quadratic extensions.
Then D. Goldfeld [G2] went quite further by employing L-functions
of elliptic curves.  These L-functions are suspected to have
central zero of order equal to the rank of the group of rational
points on the curve (the Birch and Swinnerton-Dyer conjecture).
Subsequently B. Gross and D. Zagier [GZ] provided an elliptic
curve of analytic rank three which completed Goldfeld's work with
the following estimate $$ h \gg \prod_{p|q} \, \bigl(
1-\frac{2}{\sqrt p} \bigr) \,\log q. \tag 1.15 $$ This is the
first and so far the only unconditional estimate (apart from the
implied constant, see [O]) which shows that $h \to \infty$
effectively.  Recently P. Sarnak and A. Zaharescu [SZ] used the
same elliptic curve to show that $h \gg q^{1/10}$ with an
effective constant.  However their result is conditional; they
assume (among a few minor restrictions on $q$) that $Lan_E(s) =
L_E(s) \,L_E (s,\chi)$ has no complex zeros off the critical line,
whereas the real zeros can be anywhere.
\medskip
After having exploited the power of the  central zero it seems
promising to focus on the critical line and ask if some clustering
of zeros has any effect on the class number? In fact this
possibility was contemplated in the literature independently of
the central zero effects.    In this paper we derive quite strong
and effective lower bounds for $h$, though  conditionally subject
to the existence of many small (subnormal) gaps  between zeros of
the L-function associated with a class group character. Let $$
L(s,\psi) = \sum_{\frak a} \psi({\frak a})(N{\frak a})^{-s} =
\sum_n \lambda(n)n^{-s}  \tag 1.16 $$ for $\psi \in {\widehat{\Cal
C}{\ell}}(K)$, where ${\frak a}$ runs over the non-zero integral
ideals.  This Hecke $L$-function does not factor as the Landau
product (1.11) (unless $\psi$ is a genus character), yet the
crucial feature - the lacunarity of the coefficients $$ \lambda(n)
= \sum_{N{\frak a} =n}  \psi ({\frak a})  \tag 1.17 $$ - appears
if the class number is assumed to be relatively small. One can
show that the number of zeros of $L(s,\psi)$ in the rectangle
$s=\sigma + it$ with $0 \leqslant \sigma \leqslant 1, \,0 < t
\leqslant T$ satisfies $$ N(T,\psi) = \frac{T}{\pi} \,\log
\,\frac{T\sqrt q}{2\pi} - \frac{T}{\pi} + O(\log qT).  \tag 1.18
$$ This indicates (assuming GRH for $L(s,\psi)$) that the average
gap between consecutive zeros $\rho = \fc + i \gamma$ and
$\rho^\prime = \fc + i\gamma^\prime$ is about $\pi/\log \gamma$.
\medskip
We prove that if the gap is somewhat smaller than the average for
sufficiently many pairs of zeros on the critical line (no Riemann
hypothesis is required) then $h \gg \sqrt q (\log q)^{-A}$ for
some constant $A > 0$. Actually we establish various more general
results among which are the following two theorems.  Let $\rho =
\fc + i\gamma$ denote the zeros of $L(s,\psi)$ on the critical
line and $\rho^\prime = \fc + i\gamma^\prime$ denote the nearest
zero to $\rho$ on the critical line (we assume that $\rho^\prime
\neq \rho$ except when $\rho$ is a multiple zero in which case
$\rho' = \rho$).  Note that we do not count zeros off the critical
line, but we allow them to exist. For  $0 < \alpha \leqslant 1 $
and $T \geqslant 2 $ we put $$ D(\alpha,T) = \#\{\rho;\,\,\, 2
\leqslant \gamma \leqslant T,\,\, \,|\gamma - \gamma^\prime|
\leqslant \frac{\pi(1-\alpha)}{\log \gamma} \} \tag 1.19 $$
\proclaim{Theorem 1.1}  Let $A \geqslant 0$ and $\log T \geqslant
(\log q)^{A+6}$. Suppose $$ D(\alpha,T) \geqslant \frac{c T \,\log
\,T}{\alpha\,(\log q)^A}  \tag 1.20 $$ for some $0 <\alpha
\leqslant 1$, where $c$ is a large absolute constant. Then
\endproclaim $$ L(1,\chi) \geqslant (\log T)^{-2} (\log q)^{-2A
-6}.  \tag 1.21 $$
\medskip
This result is a special case of Proposition 10.1.  Taking $A=12$ and
$\log T = (\log q)^{18}$ we get $L(1,\chi) \geqslant (\log q)^{-66}$,
provided $D(\alpha,T) \geqslant \alpha^{-1} c \,T (\log \,T)^\frac 1 3$.
\medskip
An interesting case is $\zeta_K(s) = \zeta(s) L(s,\chi)$ (where
$\zeta(s)$ is the Riemann zeta function), that is the case of the
trivial class group character.  Since we do not need all the zeros
we choose only these of $\zeta(s)$ and state the conditions in
absolute terms (without mentioning the conductor $q$, see
Corollary 10.2).
\medskip
\proclaim{Theorem 1.2}  Let $\rho = \fc + i\gamma$ be the zeros of
$\zeta(s)$ on the critical line and $\rho^\prime = \fc +
i\gamma^\prime$ be the nearest zero to $\rho$ on the critical line
($\rho' = \rho$ if $\rho$ is multiple ). Suppose $$ \# \{\rho;\,\,
\,0 < \gamma \leqslant T,\,\, |\gamma - \gamma^\prime| \leqslant
\frac{\pi}{\log\,\gamma}(1- \frac{1}{\sqrt {\log \, \gamma}})\}
\gg T(\log \,T)^\frac 4 5 \tag 1.22 $$ for any $T \geqslant 2001$.
Then we have $$ L(1,\chi) \gg (\log q)^{-90} \tag 1.23 $$ where
the implied constant is effectively computable.
\endproclaim
\medskip
  Many other results can be inferred from Proposition 10.1. We
selected our points $\rho,\rho^\prime$ from zeros on the critical
line.  However it is not hard to include other zeros in the
critical strip, or even points where $L(s,\psi)$ or even
$L^\prime(s,\psi)$ is small. As an illustration, the following
assertion follows immediately from Proposition 10.1.
\medskip
\proclaim{Corollary 1.3}  Suppose there are points  $2 \leqslant
t_1 < \ldots < t_R \leqslant T $ with  $t_{r+1} - t_r \geqslant 1
$ such that$$|L^\prime({\tfrac 1 2} +it_r,\psi)| \leqslant (\log
q)^\frac 7 2 $$ for $r=1,...,R$,  where $R \gg T= {\exp}(\log\,
q)^6$. Then $$L(1,\chi) \gg (\log q)^{-18}.$$
\endproclaim
\medskip
Considerations of Random Matrix Theory (see [Hu]) suggest that the
hypothesis above is likely to be achieved.
\medskip

 Many sections of this
paper are valid for arbitrary points in the strip (not necessarily
zeros of $L(s,\psi)$); it is only in the last four sections that
we select the points on the line Re $s = \fc$ to simplify the
arguments. Thus, if (1.22) were established unconditionally for
pairs of zeros $\rho = \beta + i \gamma,\,\,\rho' = \beta' + i
\gamma'$, which may or may not be on the critical line, then
(1.23) would hold. On the other hand one should be careful of
charging the Riemann hypothesis. Although Theorem 1.2 does not
require the Riemann hypothesis, we can imagine that someone shows
the condition (1.22) using the Riemann hypothesis. In this
scenario one still cannot conclude unconditional, effective bound
(1.23).
\medskip
Note that the average gap between consecutive zeros of $\zeta(s)$
is $2\pi/\log \gamma$, so we count in (1.22) the gaps which are
slightly smaller than the half of the average.  In view of the
implications for the class number, one has a good reason to search
for small gaps between zeros of $\zeta(s)$.  This task was
undertaken long ago.  Let us assume the Riemann hypothesis for
$\zeta(s)$. First H.L. Montgomery [M] showed that $$ |\gamma -
\gamma^\prime| < \frac{2\pi\theta}{\log \gamma}  \tag 1.24 $$
infinitely often with $\theta = 0.68$.  This was subsequently
lowered to $\theta = 0.5179$ by Montgomery and Odlyzko [MO], to
$\theta = 0.5171$ by Conrey, Ghosh and Gonek [CGG], and to $\theta
= 0.5169$ by Conrey and Iwaniec (work in progress).  We doubt that
the current technology is capable to reduce (1.24) down to $\theta
= \fc$.  Nevertheless it is an attractive and realistic
proposition to get (1.24) with any $\theta > \fc$.
\medskip
The well justified Pair Correlation Conjecture (PCC) of H.
Montgomery [M] does imply (1.24) with any $\theta > 0$ for a
positive density of zeros.  Precisely one expects that $$ \align
\# \bigl\{ m \ne n;\quad &0 < \gamma_m,\gamma_n \leqslant T,\quad
\frac{2 \pi \alpha}{\log \,T} < \gamma_m - \gamma_n < \frac{2 \pi
\beta}{\log \,T} \bigr\} \\
  &\sim \frac{T}{2 \pi} (\log \,T) \int^\beta_\alpha \bigl( 1 -
(\frac{\text{sin }\pi u}{\pi u})^2 \bigr) \, du
\endalign
$$ as $T \rightarrow \infty$, for any fixed $\beta > \alpha $.
\medskip
\remark{Remarks}  Montgomery says he was led to formulate the PCC when
looking for small gaps between zeros of $\zeta(s)$ in connection to
the class number problem.  We cannot guess how precise the connection
he established at that time.  However, Montgomery did publish a joint
paper with P.J. Weinberger [MW] in which they used zeros of fixed real
$L$-functions close to the central point $s = \frac{1}{2}$ to derive
explicit estimates and to perform extensive numerical computations for the
imaginary quadratic fields $K = \Bbb Q(\sqrt{-q})$ with the class
number $h = 1,2$.

In a similar fashion M. Jutila [J] considered a large family of
Dirichlet L-functions $L(s,\chi)$ for all $\chi$(mod $k$) with $k
\leqslant X$, and he showed that their zeros near the central
point tend to form an arithmetic progression if the class number
of $K = \Bbb Q(\sqrt -q)$ is relatively small.  Our principal idea
in this paper is reminiscent of that of Jutila .  However, as we
are dealing with a single L-function (no averaging over
characters) our arguments are quite intricate, especially when we
have to deal with the off-diagonal terms in the mean-value of
$|L(\fc + it, \psi)|^2$ (see Theorem 6.1).  Jutila's arguments do
not go that far.

Our general result in Proposition 10.1 would also imply approximate
periodicity in the distribution of most of the zeros of $\zeta(s)$
(which is inherited from oscillation of the root number (7.26)),
if we have assumed that the class number was small.  This clearly
violates the distribution law of zeros according to the PCC.
\medskip
At the meeting in Seattle  of August 1996 R. Heath-Brown gave a
lecture ``Small Class Number and the Pair Correlation of Zeros''
in which he communicated results (still unpublished) some of which
are similar to ours, yet they are more restrictive.  Heath-Brown
requires $L(1,\chi) \ll q^{-\delta}$ for some constant $\delta
> \frac{1}{4}$, which condition contradicts the Siegel bound
$L(1,\chi) \gg q^{-\var}$, but his arguments are effective so the
results remain valid .
\endremark
\mn {\it Acknowledgement.}  This work began during the second
author visit to the American Institute of Mathematics in summer
1999. He has pleasure to acknowledge support and the hospitality
of the Institute.  The final version was written during the second
author visit to the University of Lille in June 2001, and he is
thankful for this opportunity. \mn {\it Note added in May 2001.}
We found our results in Section 3 and Section 4 to be similar to
these in Appendix A and Appendix B of the paper ``Rankin-Selberg
$L$-functions in the level aspect" by E. Kowalski, P. Michel and
J. Vanderkam (to appear).  Had we known their results earlier we
would gladly incorporate them to reduce our arguments. However, we
decided not to modify our original parts to preserve the
self-contained presentation.
\medskip
\subhead 2.  Basic Automorphic Forms\endsubhead
\medskip
We are mainly interested in L-functions for characters on ideals
in the imaginary quadratic field $K = \Bbb Q (\sqrt -q)$.  Every
such L-function is associated with a holomorphic automorphic form
of level $q$ and the real primitive character $\chi$(mod $q$) (the
nebentypus). However, to get better perspective, we begin by
reviewing the whole spectrum  of real-analytic forms.  In
particular we focus on the Eisenstein series, because they are
most important automorphic forms for our applications to Dirichlet
L-functions (they correspond to genus characters of the class
group of $K$).  Some more details and proofs can be found in [I]
and [DFI2].

The group $SL_2(\Bbb R)$ acts on the upper-half plane $\Bbb H$ by
the linear fractional transformations $\gamma z = (az + b)/(cz +
d)$ if $\gamma = \left (\smallmatrix a&b\\c&d\endsmallmatrix
\right) \in SL_2(\Bbb R)$.  We put $$ j_\gamma(z) = \frac{cz +
d}{|cz + d|}. $$ Note that $j_{\beta \gamma}(z) = j_\beta (\gamma
z) j_\gamma(z)$.  Next we fix a positive integer $k$ and put $$
J_\gamma(z,s) = j^{-k}_\gamma(z) (\text{Im } \,\gamma z)^s  =
\bigl( \frac{cz +d}{|cz +d|}\bigr)^{-k} \bigl(\frac{y}{|cz +
d|^2}\bigr)^s $$ for $\gamma \in SL_2 (\Bbb R), z \in \Bbb H$ and
$s \in \Bbb C$.  Note that $J_{\beta \gamma} (z,s) = j^{-k}_\gamma
(z) \,J_\beta (z,s)$. Since $ J_\gamma (z,s) $ depends only on the
lower row $(c,d)$ of $\gamma$ we shall write $J_{(c,d)} (z,s)$ in
place of $J_\gamma (z,s)$. Actually $J_{(c,d)}(z,s)$ is defined by
the last expression for any pair of real numers $c,d,$ not both
zero. Note that for $u > 0$ we have $J_{(uc,ud)} (z,s) = u^{-2s}
J_{(c,d)} (z,s)$.
\medskip
Throughout $\Gamma = \Gamma_0 (q)$ denotes the Hecke congruence
group of level $q$; its index in the modular group is $$ \nu(q) =
[\Gamma_0 (1): \Gamma_0 (q)] = q \prod_{p|q} (1 + \frac 1 p). \tag
2.1 $$ To simplify the presentation (without compromising our
applications) we restrict $q$ to odd, squarefree numbers.  Let
$\chi = \chi_q$ be the real primitive character of conductor $q$,
i.e. $\chi_q (n) = (\frac n q)$ is the Jacobi-Legendre symbol.
This induces a character on $\Gamma$ by $$ \chi (\gamma) = \chi
(d), \,\,\,\,\,\,\, \text{ if } \gamma = \left (\smallmatrix
a&b\\c&d
\endsmallmatrix \right ) \in \Gamma.  \tag 2.2 $$
\medskip
We are interested in the space ${\Cal A}_k(\Gamma,\chi)$ of
automorphic functions of weight $k \geqslant 1$ for the group
$\Gamma $ and character $\chi $, i.e. the functions $f:\Bbb H \to
\Bbb C$ satisfying $$ f(\gamma z) = \chi(\gamma) j ^k_\gamma(z)
f(z), \quad \text{ if } \gamma \in \Gamma.  \tag 2.3 $$ We assume
$\chi (-1) = (-1)^k$, or otherwise ${\Cal A}_k (\Gamma, \chi)$
consists only of the zero function. The Laplace operator $$
\Delta_k = y^2 (\frac{\partial^2}{\partial x^2} +
\frac{\partial^2}{\partial y^2}) - iky \frac{\partial}{\partial x}
$$ acts on ${\Cal A}^\infty_k (\Gamma,\chi)$ - the subspace of
smooth automorphic functions. Any $f \in {\Cal A}^\infty_k
(\Gamma,\chi)$ which is eigenfunction of $\Delta_k$, say
$(\Delta_k + \lambda)f = 0$, is called Maass form of eigenvalue
$\lambda$.

Our primary examples of Maass forms are the Eisenstein series
associated with cusps of $\Gamma$.  Let $\Gamma_{\frak a} =
\{\gamma \in \Gamma;\, \gamma {\frak a} = {\frak a}\}$ be the
stability group of the cusp ${\frak a}$.  There exists
$\sigma_{\frak a} \in SL_2(\Bbb R)$ such that $\sigma_{\frak a}
\infty = {\frak a}$ and $\sigma^{-1}_{\frak a}\, \Gamma_{\frak
a}\, \sigma_{\frak a}\, = \Gamma_\infty$ - the group of
translations by integers.  We call $\sigma_{\frak a}$ a scaling
matrix of ${\frak a}$. The Eisenstein series associated with
${\frak a}$ is defined by $$ E_{\frak a}(z,s) = \sum_{\gamma \in
\Gamma_{\frak a} /\Gamma} \chi (\gamma) J_{\sigma^{-1}_{\frak a}
\gamma} (z,s).  \tag 2.4 $$ This series converges absolutely for
Re $s>1$, it does not depend on the choice of $\sigma_{\frak a}$,
nor on the choice of ${\frak a}$ in its equivalence class.  The
Eisenstein series $E_{\frak a} (z,s)$ is a Maass form of
eigenvalue $\lambda = s (1-s)$.

Any cusp ${\frak a}$ of $\Gamma = \Gamma_0 (q)$ is equivalent to a rational
point $1/v$, where $v$ is a divisor of $q$ (recall that $q$ is
squarefree).  Put
$$
q = uw  \tag 2.5
$$
so $w$ is the width of the cusp ${\frak a} \sim 1/v$.  As a scaling matrix of
${\frak a} \sim 1/v$ we can choose
$$
\sigma_{\frak a} = \left ( \matrix \sqrt w & 0\\ v \sqrt w & 1/ \sqrt w
\endmatrix \right ).
$$ Next, according to (2.5), we factor the character $\chi_q =
\chi_v \chi_w$.  Then the Eisenstein series (2.4) can be written
explicitly as follows $$ E_{\frak a} (z,s) = \frac{1}{2w^s}
\underset{(c,d) = 1}\to {\dsize\sum\sum} \chi_v (d) \chi_w (-c)
J_{(cv,d)} (z,s)  \tag 2.6 $$ where $c,d$ run over co-prime
integers.  Hence applying Poisson's summation one can derive a
Fourier expansion of $E_{\frak a} (z,s)$ (in terms of the
Whittaker function) from which one can see (among other things)
that $E_{\frak a} (z,s)$ is meromorphic in the whole complex
s-plane without poles in Re $s \geqslant \fc$ (see (7.12) and
(7.13) of [DFI2]).

The Eisenstein series $E_{\frak a}(z,s)$ on the line Re $s=\fc$
yield an eigenpacket of the continous spectrum of $\Delta_k$ in
the subspace ${\Cal L}_k (\Gamma,\chi)$ of square-integrable
functions $f(z) \in {\Cal A}_k (\Gamma,\chi)$ with respect to the
invariant measure $y^{-2} dxdy$.  The continous spectrum covers
the segment $[\frac 1 4, \infty)$ with multiplicity $\tau (q)$
(the number of inequivalent cusps equals the number of divisors of
$q$).  Let $\Cal E_k (\Gamma,\chi) \subset {\Cal L}_k
(\Gamma,\chi)$ be the subspace of the continous spectrum (it is a
linear space spanned by a kind of incomplete Eisenstein series).
Let $\Cal C_k (\Gamma,\chi)$ be the orthogonal complement of $\Cal
E_k (\Gamma,\chi)$ in $\Cal L_k (\Gamma,\chi)$, so $\Cal
L_k(\Gamma,\chi) = \Cal E_k (\Gamma,\chi) \oplus \Cal C_k
(\Gamma,\chi)$. The Laplace operator $\Delta_k$ acts on $\Cal C_k
(\Gamma,\chi)$, and it has an infinite, purely discrete spectrum
in the segment $[\frac k 2 (1- \frac k 2), \infty)$.  In other
words $\Cal C_k (\Gamma,\chi)$ is spanned by square-integrable
automorphic forms.  These are characterized by vanishing at every
cusp (because they are orthogonal to every incomplete Eisenstein
series), and are called Maass cusp forms.

>From now on we take only the Maass cusp forms $f(z)$ of the Laplace
eigenvalue $\lambda = \frac k 2(1- \frac k 2)$, and if $k = 1$ we also
take the Eisenstein series $E_{\frak a}(z,s)$ at $s = \fc$.  All these forms
come from the classical holomorphic forms of weight $k$, precisely we
have
$$
F(z) = y^{-\frac k 2} f (z) \in S_k (\Gamma,\chi)
$$
$$
E_{\frak a}(z) = y^{-\fc} E_{\frak a} (z,\tfrac 1 2) \in M_1(\Gamma,\chi).
$$
For any $n \geqslant 1$ the Hecke operator $T_n$ is defined on
$M_k(\Gamma,\chi)$ by
$$
(T_n F)(z) = \frac{1}{\sqrt n} \sum_{ad = n} \chi(a) (\frac a d)^{k/2}
\sum_{b(\text{mod } d)} F(\frac{az + b}{d}).
$$
There is a basis of $S_k(\Gamma,\chi)$ which consists of eigenforms of
the Hecke operators $T_n$ with $(n,q) = 1$.  Moreover, by the
multiplicity-one property (which holds in $S_k(\Gamma,\chi)$ because
$\chi$ is primitive of conductor equal to the level) these forms are
automatically eigenfunctions of all $T_n$.  Consequently we may assume
that
$$
T_n F = \lambda_F (n) F \qquad \text{ for all } n \geqslant 1.  \tag 2.7
$$
After a normalization of $F$ by a suitable scalar the Hecke
eigenvalues $\lambda_F(n)$ agree with the coefficients in the Fourier
series
$$
F(z) = \sum^\infty_1 \lambda_F (n) n^{\frac{k-1}{2}} e (nz).  \tag 2.8
$$
Such $F$ is called a primitive cusp form of weight $k$, level $q$ and
character $\chi$.
\medskip
One can show that the modified Eisenstein series $y^{-\frac k 2}
E_{\frak a}(z,s)$ are also eigenfunctions of all the Hecke
operators $T_n$ (see Section 6 of [DFI2]), but we are only
interested in $E_{\frak a}(z) = y^{-\fc} E_{\frak a}(z,\fc)$.  In
this case $(k=1, s=\fc)$ we have $$ T_n E_{\frak a} =
\lambda_{\frak a}(n) E_{\frak a} \qquad \text{ for all } n
\geqslant 1  \tag 2.9 $$ with $$ \lambda_{\frak a}(n) =
\sum_{n_1n_2=n} \chi_v (n_1) \chi_w (n_2).  \tag 2.10 $$ Moreover
the Hecke eigenvalues $\lambda_{\frak a}(n)$ are proportional to
the Fourier coefficients of $E_{\frak a}(z)$, specifically we have
(see [I] and [DFI2]) $$ E_{\frak a}(z) = \bar \var_v \frac{2i}{h}
\sum^\infty_0 \lambda_{\frak a}(n) \,e \,(nz) \tag 2.11 $$ where
$\var_v = \tau (\chi_v)/\sqrt v$,\, so $\var_v = 1$, $i$ according
to $v \equiv 1,3$(mod $4$) and $h$ is the class number of $K =
\Bbb Q(\sqrt {-q})$, $$ h = \pi^{-1} \sqrt q \, L(1,\chi).  \tag
2.12 $$ The zero coefficient is given by $$ \lambda_{\frak a}(0)
=\cases h/2 &\text{ if } {\frak a} \sim \infty, 0\\ 0
&\text{otherwise }\endcases.  \tag 2.13 $$ Our particular
Eisenstein series $E_{\frak a}(z)$ (recall that in this case we
have $k = 1$ and $\chi_q (-1) = -1$ so $ q \equiv 3$(mod $4$)) can
be expressed by theta functions for ideal classes of $K = \Bbb Q
(\sqrt {-q})$. For every class ${\Cal A} \in {\Cal C}\ell(K)$ we
put $$ \theta_{\Cal A} (z) = \fc + \sum_{{\frak a} \in {\Cal A}} e
(z N {\frak a}) \tag 2.14 $$ where ${\frak a}$ runs over integral
ideals in ${\Cal A}$ and $N {\frak a}$ is the norm of ${\frak a}$
(the number 2 stands for the number of units, we assume $q \neq
3)$. This theta function is also given by $$ \theta_{\Cal A}(z) =
\fc \sum_m \sum_n e (z \varphi_{\Cal A}(m,n))  \tag 2.15 $$ where
$\varphi_{\Cal A} (x,y) = ax^2 + bxy +cy^2$ is the corresponding
quadratic form.  Specifically we have $a > 0, (a,b,c) = 1, \,b^2 -
4 ac = -q$ and $$ {\frak a} = a \,\Bbb Z + \frac{b + i \sqrt q}{2}
\,\Bbb Z $$ is an integral primitive ideal representing the class
${\Cal A}$.  One shows that the theta function $\theta_{\Cal
A}(z)$ for any class ${\Cal A}$ belongs to $M_1(\Gamma,\chi)$.
Hence for any character $\psi \in {\widehat{\Cal C}\ell}(K)$ $$
\theta(z;\psi) = \sum_{{\Cal A} \in {\Cal C}\ell(K)} \psi({\Cal
A}) \theta_{\Cal A}(z) \tag 2.16 $$ is an automorphic form of
weight one, level $q$ and character $\chi = \chi_q$.  Note that
$\theta (z;\bar \psi) = \theta (z;\psi)$.  This has the Fourier
expansion $$ \theta(z;\psi) = \sum^\infty_0 \lambda_\psi (n) \,e
\,(nz)  \tag 2.17 $$ with $\lambda_\psi (0) = \delta_\psi h/2$,
and for $n \geqslant 1$ $$ \lambda_\psi (n) = \sum_{N {\frak a} =
n} \psi({\frak a}).  \tag 2.18 $$ In particular the Eisenstein
series $E_{\frak a}(z)$ are obtained from theta functions for real
class group characters.  Any real character $\psi \in
{\widehat{\Cal C}\ell}(K)$ is given uniquely by $$ \psi({\frak p})
= \cases \chi_v(N{\frak p}), \quad \text{ if } p \nmid v \\
\chi_w(N{\frak p}),\quad \text{ if } p \nmid w \endcases \tag 2.19
$$ where $\chi_v \chi_w = \chi_q$ (note that $\psi({\frak a})$ is
well defined by (2.19) because $\chi_q (N{\frak a}) = 1$ if
$({\frak a},q) = 1$).  Interchanging $v$ and $w$ we obtain the
same $\psi$.  However different factorizations $vw = q$ up to the
order yield distinct real class group characters. Therefore we
have exactly $\fc \tau(q)$ such characters, they are called the
genus characters.  If ${\frak a} \sim 1/v$ then we have $$
E_{\frak a}(z) = { {\bar{\var}}_v} \frac{2i}{h} \theta (z;\psi)
\tag 2.20 $$ where $\psi \in {\widehat{\Cal C}\ell} (K)$ is the
genus character given by (2.19) and $\lambda_{\frak a}(n) =
\lambda_\psi(n)$ for all $n \geqslant 0$ (see (2.11)). Note that
the Eisentein series $E_{\frak a}(z)$ and $E_{{\frak a}'}(z)$ for
the ``transposed'' cusps ${\frak a} \sim 1/v$ and ${\frak
a}^\prime \sim 1/w$ are linearly dependent, in fact $\var_v
\,E_{\frak a}(z) = \var_w \,E_{{\frak a}^\prime}(z)$ (this is true
only for the Eisenstein series at the central point $s = \fc$!).

If $\psi \in {\widehat{\Cal C}\ell}(K)$ is not real, then the theta function
$\theta (z;\psi)$ is a primitive cusp form of weight one with Hecke
eigenvalues $\lambda_\psi (n)$ given by (2.18).
\medskip
Cusp forms of any odd weight can be constructed from the class group
characters as follows.  Let $k$ be odd, $k > 1$ and $q \equiv 3$(mod
$4$), $q > 3$.  Let $\psi$ be a character on ideals in $K = \Bbb Q(\sqrt
{-q})$ such that
$$
\psi ((\alpha)) = \biggl( \frac{\alpha}{|\alpha|} \biggr)^{k-1}  \tag 2.21
$$
for any $\alpha \in K^{\ast}$.  All such characters are obtained by
multiplying a fixed character with the class group characters so we
have exactly $h = h(-q)$ characters of type (2.21) (we say of frequency
$k-1$).  For every $\psi$ of frequency $k-1$ we associate the function
$$
\theta (z;\psi) = \sum_{\frak a} \psi ({\frak a}) (N{\frak
a})^{\frac{k-1}{2}} \,e \,(zN{\frak a})  \tag
2.22
$$
where ${\frak a}$ runs over the non-zero integral ideals.  One shows that
$\theta(z;\psi) \in S_k(\Gamma,\chi)$ and that $\theta(z;\psi)$ is a
primitive cusp form with Hecke eigenvalues $\lambda_\psi (n)$ given by
(2.18)(see Section 12.3 of [I]).
\medskip
Besides (2.3) and (2.7) the primitive forms satisfy some bilateral
modular equations which are obtained by certain transformations $\omega
\in SL_2 (\Bbb R)$ not in the group $\Gamma_0(q)$. For any $\omega =
\left ( \smallmatrix a&b\\c&d \endsmallmatrix \right ) \in SL_2 (\Bbb
R)$ the $\omega$-stroke operator is defined on functions $F:\Bbb H \to
\Bbb C$ by
$$
F_{|\omega}(z) = (cz + d)^{-k} F(z).  \tag 2.23
$$
Note that $(F_{|\tau})_{|\sigma} = F_{|\tau \sigma}$ for any
$\tau, \sigma \in SL_2 (\Bbb R)$.
\medskip
Let $q = rs$ (recall that $q$ is squarefree so $(r,s) = 1$).  We are
interested in the $\omega$-stroke operator for
$$
\omega = \left (\matrix \alpha \sqrt r & \beta/\sqrt r\\ \gamma s
\sqrt r & \delta \sqrt r \endmatrix \right )  \tag 2.24
$$
with $\alpha, \beta, \gamma, \delta$ integers such that $\det \,\omega =
\alpha \delta r - \beta \gamma s = 1$.
\medskip
First for the $\omega$ given by (2.24) one checks that the
$\omega$-stroke maps $S_k(\Gamma,\chi)$ to itself.  Next note that
the $\omega$-stroke on $S_k(\Gamma,\chi)$ is a pseudo-involution,
precisely $$ F_{|\omega^2} = \chi_r(-1) \chi_s(r) F \tag 2.25 $$
where $\chi_r \chi_s = \chi_q$.  Moreover the $\omega$-stroke on
$S_k (\Gamma,\chi)$ almost commutes with the Hecke operators $T_n$
for $(n,q) = 1$, precisely $$ T_n (F_{|\omega}) = \chi_r(n) (T_n
F)_{|\omega}, \quad \text{ if } (n,q) = 1.  \tag 2.26 $$ Hence it
follows that if $F$ is a Hecke form (i.e. $F$ is an eigenfunction
of every $T_n$ with $(n,q)=1)$, then so is $F_{|\omega}$ (of
course, with different Hecke eigenvalues).  By the
multiplicity-one property it follows that both $F$ and
$F_{|\omega}$ are primitive (i.e. the eigenfunctions of all
$T_n$).  Therefore for any primitive form $F \in S_k
(\Gamma,\chi)$ there exists a unique primitive form $G \in S_k
(\Gamma,\chi)$ and a complex number $\eta_F(\omega)$ such that $$
F_{|_\omega} = \eta_F (\omega) G.  \tag 2.27 $$ As in [AL] we call
$\eta_F (\omega)$ the pseudo-eigenvalue of $|_\omega$ at $F$.  By
(2.25) we find that $\eta_F \eta_G = \chi_r (-1) \chi_s (r) = \pm
1$.  One can show that the Hecke eigenvalues of $F$ and $G$
satisfy $$ \align \lambda_G (n) = \chi_r (n) \lambda_F (n), \quad
\text{ if } (n,r) = 1, \tag 2.28 \\ \lambda_G (n) = \chi_s (n)
\bar \lambda_F (n), \quad \text{ if } (n,s) = 1.  \tag 2.29
\endalign
$$
These formulas are consistent by the property $\bar \lambda_F (n) =
\chi (n) \,\lambda_F (n)$ if $(n,q) = 1$, and they determine $G$ in terms of
$F$.  In particular we have $|\lambda_F (n)| = |\lambda_G (n)|$ for
all $n \geqslant 1$.
Hence one derives that $<G,G> \,= \,<F,F>$ \,and \,$<F,F> \,=
\,|\eta_F(\omega)|^2 <G,G>$, so
$$
|\eta_F (\omega)| = 1.  \tag 2.30
$$
\medskip
Note that $G$ depends only on $r,s$ ($G$ is a hybrid twist of $F$
by the characters $\chi_r, \chi_s$), but not on $\alpha, \beta,
\gamma, \delta$ in $\omega$.  If $\omega$ and $\omega^\prime$ are
given by (2.23) with the same $r,s$ then $$ \rho = \omega'
\omega^{-1} = \left (\matrix \alpha^\prime \delta r- \beta^\prime
\gamma s & \beta^\prime \alpha - \alpha^\prime \beta\\
(\gamma^\prime \delta - \delta^\prime \gamma) q & \delta^\prime
\alpha r - \gamma^\prime \beta s \endmatrix \right ) $$ $$
F_{|_{\omega^{\prime}}} = F_{|\rho \omega} = (F_{|\rho})_{|\omega}
= \chi(\rho) \,F_{|\omega} = \chi (\rho) \eta_F (\omega) G $$ $$
\chi(\rho) = \chi (\delta^\prime \alpha r - \gamma^\prime \beta s)
= \chi_r (-\gamma^\prime \beta s) \chi_s (\delta^\prime \alpha r)
= \chi_r (\beta^\prime/\beta) \chi_s (\alpha^\prime/\alpha) $$ by
the determinant equation $\alpha^\prime \delta^\prime r -
\beta^\prime \gamma^\prime s = 1$. Hence we get the relation
$\eta_F (\omega^\prime) = \chi_r (\beta^\prime \beta) \chi_s
(\alpha^\prime \alpha) \eta_F (\omega)$.  This relation shows that
the pseudo-eigenvalue $\eta_F (\omega)$ of $\omega$ given by
(2.23) factors into $$ \eta_F (\omega) = \chi_r (\beta) \chi_s
(\alpha) \eta_F (r,s) = \chi_r (-\gamma s) \chi_s (\delta r)
\eta_F (r,s)  \tag 2.31 $$ where $\eta_F(r,s)$ depends only on
$r,s$ and $F$.
\medskip
The case $r = q$ and $s = 1$ is special.  We can choose
$$
\omega = \left (\smallmatrix 0 & -1/\sqrt q\\ \sqrt q & 0
\endsmallmatrix \right )  \tag 2.32
$$
getting $F_{|_\omega} = \eta_F \,\bar F$, where $\bar F$ is obtained
from $F$ by complex conjugating the coefficients in the Fourier
expansion (2.8).  Moreover in this case one shows that (see Theorem
6.29 of [I])
$$
\eta_F = \var_q \bar \lambda_F (q).  \tag 2.33
$$
\medskip
The modified Eisenstein series $y^{-\frac k 2} E_{\frak a}(z,s)$ is also a
pseudo-eigenfunction of the $\omega$-stroke operator.  We shall verify
this fact by explicit computations rather than by going through the
theory of Hecke operators.
Although we are only interested in $E_{\frak a}(z)
= y^{-\fc} E_{\frak a}(z,\fc)$ we present the computations in a
general case (i.e. for any $k \geqslant 1$) for record.
Note that for any $\gamma, \omega \in SL_2 (\Bbb R)$
$$
y^{-\frac k 2} J_\gamma (z,s)_{|_\omega} = y^{-\frac k 2} J_{\gamma
\omega} (z,s).
$$
Hence we get by (2.6) for Re $s > 1$
$$
y^{-\frac k 2} E_{\frak a}(z,s)_{|_\omega} = \tfrac{1}{2}
w^{-s} y ^{-\frac k 2}
\underset{(c,d) = 1}\to {\dsize\sum\sum} \chi_r (d) \chi_w (-c) J_\tau (z,s)
$$
where
$$
\tau = \left (\matrix \ast & \ast\\ cv & d \endmatrix \right
) \omega = \left (\matrix \ast & \ast\\ (\alpha c v + \gamma d s)
\sqrt r & (\beta c v + \delta d r) /\sqrt r \endmatrix \right ).
$$
Put $r^\prime = ( r,v), \,r_1 = r/(r,v), \,r_2 = v/(r,v)$ and $s^\prime =
(s,v),s_1 = s/(s,v), \,s_2 = v/(s,v)$.  Since $q = vw = rs$ is
squarefree we have $r_1s_1 = w, r_2s_2 = v, r_1s_2 = r, r_2s_1 = s$.
In the lower row of $\tau$ we extract the factors $r^\prime, s^\prime$
getting
$$
\tau = \left (\matrix \ast&\ast\\
C s^\prime \sqrt r \,& \,D r^\prime /\sqrt r \endmatrix \right )
$$
where $C = \alpha c s_2 + \gamma d s_1$ and $D = \beta c r_2 + \delta
d r_1$.  Solving this system of linear equations of determinant
$\alpha \delta s_2 r_1 - \beta \gamma r_2 s_1 = \alpha \delta r -
\beta \gamma s = 1$ we find $- c = \gamma s_1 \,D - \delta r_1 \,C$ and $d
= \alpha s_2 D - \beta r_2 C$.  Hence the condition $(c,d) = 1$ is
equivalent to $(C,D) = 1$.  Next we factor the characters $\chi_v =
\chi_{r_2} \chi_{s_2}$ and $\chi_w = \chi_{r_1} \chi_{s_1}$ to compute
$\chi_v (d) = \chi_{r_2} (\alpha s_2 D) \chi_{s_2} (-\beta r_2 C)$ and
$\chi_w (-c) = \chi_{r_1} (\gamma s_1 \,D) \chi_{s_1} (-\delta r_1 C)$.
Hence $\chi_v (d) \chi_w (-c) = \chi_{r_1 r_2} (D) \chi_{s_1 s_2} (-C)
\,\eta$, where
$$
\eta = \chi_{r_2}(\alpha s_2) \chi_{s_2}(\beta r_2) \chi_{r_1}(\gamma
s_1) \chi_{s_1}(\delta r_1).  \tag 2.34
$$
Extracting from $\tau$ the factor $u = r^\prime/\sqrt r$ by the
property $J_{(uc,ud)} (z,s) = u^{-2s} J_{(c,d)} (z,s)$, and using
$s^\prime \sqrt r/u = r_1 r_2, \,u^2w \,r_1 r_2 = q$, we conclude from
the above computations that
$$
y^{-\frac k 2} E_{\frak a}(z,s)_{|_\omega} = \eta y^{-\frac k 2}
E_{{\frak a}^\ast}(z,s)  \tag 2.35
$$
where $E_{{\frak a}^\ast}(z,s)$ is
the Eisenstein series for the cusp ${\frak a}^\ast
\sim 1/r_1 r_2$, i.e.
$$
{\frak a}^\ast \sim \frac{(r,v)^2}{r v}.  \tag 2.36
$$
By the determinant equation $\alpha \delta s_2 r_1 - \beta \gamma r_2
s_1 = 1$ we eliminate $\alpha, \beta$ in (2.34) getting
$$
\eta = \chi_r (\frac{\gamma s}{(s,v)}) \chi_s (\frac{\delta r}{(r,v)})
\chi_{v/(s,v)} (-1).  \tag 2.37
$$
In particular for $k = 1$ and $s = \fc$ we obtain from (2.35) (by
analytic continuation)
\proclaim{Proposition 2.1}  Let $q = vw = rs > 1$ be squarefree and
odd.  Then the holomorphic Eisenstein series $E_{\frak a}(z) = y^{-\fc}
E_{\frak a}(z,\fc)$ for cusp ${\frak a} \sim 1/v$ is a
pseudo-eigenfunction of the
$\omega$-stroke operator (with $\omega$ given by (2.23)), specifically
$$
E_{{\frak a}|\omega} = \eta E_{{\frak a}^\ast}  \tag 2.38
$$
where ${\frak a}^\ast$ is given by (2.36) and $\eta$ by (2.37). \endproclaim
\medskip
In the special case (2.32) the equation (2.38) becomes
$$
(z\sqrt q)^{-1} E_{1/v} (\frac{-1}{qz}) = \chi_v (-1) E_{1/w}(z). \tag
2.39
$$
\medskip
\subhead 3. Summation Formulas\endsubhead
\medskip
Suppose we have two functions $A(z), B(z)$ on $\Bbb H$ given by Fourier series
$$
\align
A(z) = \sum^\infty_0 a_n e(nz)  \tag 3.1 \\
B(z) = \sum^\infty_0 b_n e(nz)  \tag 3.2
\endalign
$$
with $a_n,\,b_n \ll n^{k-1 + \var}$.  Suppose that $A(z),\,B(z)$ are
connected by the $\omega$-stroke operator, say
$$
A_{|_\omega}\,(z) = \eta\, B(z)  \tag 3.3
$$
for some $\omega = \left (\smallmatrix a & b\\c & d \endmatrix \right
) \in SL_2 (\Bbb R)$ with $c > 0$ and some complex number $\eta \neq
0$.  In particular for $z = (-d + iy)/c$ we have
$$
(iy)^{-k} A(\frac a c + \frac{i}{cy}) = \eta\, B(\frac{-d}{c} +
\frac{iy}{c})  \tag 3.4
$$
for any $y > 0$.  Following Hecke this formula can be expressed as a
functional equation for the L-functions
$$
\align
L_A (s,\frac a c) =
\sum^\infty_1 a_n e\,(\frac{an}{c}) n^{-s}  \tag 3.5 \\
L_B (s,\frac{-d}{c}) = \sum^\infty_1 b_n e\,(\frac{-dn}{c}) n^{-s}.
\tag 3.6
\endalign
$$
Put
$$
\align
\Lambda_A (s,\frac a c) = (\frac{c}{2\pi})^s \,\Gamma (s) L_A (s,\frac
a c)  \tag 3.7 \\
\Lambda_B (s,\frac{-d}{c}) = (\frac{c}{2\pi})^s\,\Gamma (s) L_B
(s,\frac a c).  \tag 3.8
\endalign
$$
First we establish by integrating (3.4) the following formula
$$
\align
\Lambda_A(s,\frac a c) + \frac{a_0}{s} + i^k \eta \frac{b_0}{k-s} &=
\int^\infty_1 [A(\frac a c + \frac{iy}{c}) - a_0] y^{s-1} dy\\
&+ i^k \eta \int^\infty_1 [B(\frac{-d}{c} + \frac{iy}{c}) - b_0]
y^{k-s-1} dy \endalign
$$
for Re $s > k$.  Since $A(z) - a_0$ and $B(z) - b_0$ have
exponential decay as $y = \text{ Im } z \to \infty$, the above integrals
converge absolutely and they are entire functions bounded on veritical
strips.  Similarily we have (because $B_{|_{\omega^{-1}}}(z) =
\eta^{-1} A(z)$ and $\omega^{-1} = - \sssize{\left(\smallmatrix -d & \ast\\ c
& -a \endsmallmatrix\right)}$)
$$
\align
\Lambda_B(s,\frac{-d}{c}) + \frac{b_0}{s} + (i^k \eta)^{-1}
\frac{a_0}{k-s} &= \int^\infty_1 [B(\frac{-d}{c} + \frac{iy}{c}) -b_0]
y^{s-1} dy\\
&+ (i^k \eta)^{-1} \int^\infty_1 [A(\frac a c + \frac{iy}{c}) -a_0]
y^{k-s-1} dy.  \endalign
$$
Combining both formulas we obtain the following functional equation
$$
\Lambda_A(s,\frac a c ) = i^k \eta \Lambda_B (k-s,\frac{-d}{c}).  \tag
3.9
$$
\medskip
For notational convenience we put
$$
a_n = n^{\frac{k-1}{2}} a(n), \quad b_n = n^{\frac{k-1}{2}} b(n)  \tag
3.10
$$
so the corresponding L-functions are shifted from $s$ to $s-
\frac{k-1}{2}$, and the resulting functional equation connects values
at $s$ and $1-s$.
\medskip
Next we derive from (3.9) a formula for sums of type
$$
S = \sum^\infty_1 a (n) e (\frac{an}{c}) g(n)  \tag 3.10
$$
where $g(x)$ is a nice test function.
\medskip
\proclaim{Proposition 3.1}  Suppose $A(z),\, B(z)$ satisfy (3.3) for
some $\omega = \left (\smallmatrix a & \ast\\ c & d \endsmallmatrix
\right ) \in SL_2 (\Bbb R)$ with $c > 0$ and some complex number $\eta
\neq 0$.  Then for any $g(x)$ smooth and compactly supported on $\Bbb
R^+$ we have
$$
\align
\sum^\infty_1 a(n) e (\frac{an}{c}) g (n) &= 2\pi i^k \frac{\eta}{c}
\{\frac{b_0}{\Gamma(k)} \int^\infty_0 g(x) \biggl(\frac{2 \pi \sqrt
x}{c}\biggr)^{k-1} dx \tag 3.11 \\
&+ \sum^\infty_1 b(n) e (\frac{-dn}{c}) \int^\infty_0 g(x) J_{k-1}
(\frac{4\pi}{c} \sqrt {nx}) dx\} \endalign
$$
where $J_{k-1}(x)$ is the Bessel function of order $k-1$. \endproclaim
\medskip
\noindent Proof.  The left side of (3.11) is given by the contour
integral $$ S = \frac{1}{2\pi i} \int_{(\sigma)} L_A (s,\frac a c
) G(s - \frac{k-1}{2}) ds $$ where $G(s)$ denotes the Mellin
transform of $g(x)$ and $\sigma > k$. We move to Re $s = k -
\sigma$ passing a simple pole at $s = k$ with residue $i^k \eta
b_0 (2\pi/c)^k \Gamma(k)^{-1} G((k+1)/2)$ (the point $s = 0$ is
not a pole of $L_A(s,\frac a c ))$.  Then we apply the functional
equation (3.9) getting $$ \aligned S &= i^k \eta b_0
(\frac{2\pi}{c})^k \, \Gamma(k)^{-1} G(\frac{k+1}{2})\\ &+
\frac{i^k \eta}{2\pi i} \int_{(\sigma)} L_B (s,\frac{-d}{c})
(\frac{c}{2\pi})^{2s-k} \frac{\Gamma(s)}{\Gamma(k-s)}
G(\frac{k+1}{2} - s) ds.  \endaligned $$ Expanding
$L_B(s,\frac{-d}{c})$ into the Dirichlet series and integrating
termwise we get $$ \frac{1}{2\pi i} \int_{(\sigma)} =
\frac{2\pi}{c} \sum^\infty_1 b(n) e (\frac{-dn}{c}) H (\frac{2 \pi
\sqrt n}{c}) $$ where $$ H(y) = \frac {1}{2\pi i} \int_{(\sigma)}
\frac{\Gamma(s)}{\Gamma(k-s)} G (\frac{k+1}{2} - s)) y^{k-1-2s}
ds. $$ Here we can take for $\sigma$ any positive number.  If
$\sigma < \frac k 2$ we can open the Mellin transform $$ G
(\frac{k+1}{2} -s) = \int^\infty_0 g(x) x^{\frac{k-1}{2} -s} dx $$
and change the order of integration getting $$ \align H(y) &=
\int^\infty_0 g(x) \biggl(\frac{1}{2\pi i} \int_{(\sigma)}
\frac{\Gamma (s)}{\Gamma (k-s)} (\sqrt x y)^{k-1-2s} ds \biggr)
dx\\ &= \int^\infty_0 g(x) J_{k-1} (2 \sqrt x y) dx  \endalign $$
by (6.422.9) of [GR].  This yields (3.11).
\medskip
Changing variables one can write (3.11) as follows
$$
\align
\sum^\infty_1 a(n)\,e\,(\frac{an}{c}) g (\frac{2\pi n}{c}) &= i^k \eta
\{\frac{b_0}{\Gamma(k)} \int^\infty_0 g(x) (\frac{2\pi
x}{c})^{\frac{k-1}{2}} dx \tag 3.12 \\
&+ \sum^\infty_1 b(n)\,e\,(\frac{-dn}{c})\,h\,(\frac{2\pi n}{c})\}
\endalign
$$
where $h(y)$ is a Hankel-type transform
$$
h(y) = \int^\infty_0 g(x) J_{k-1} (2 \sqrt {xy}) dx.  \tag 3.13
$$
\medskip
Now we specialize Proposition 3.1 for automorphic forms.  First we
treat the cusp forms.
\medskip
\proclaim{Proposition 3.2}  Let $F \in S_k (\Gamma, \chi)$ be a Hecke
cusp form with eigenvalues $\lambda_F (n)$.  Let $c \geqslant 1$ and $(a,c) =
1$.  Then for any function $g(x)$ smooth and compactly supported on
$\Bbb R^+$ we have
$$
\align
\sum^\infty_1 \lambda_F &(n)\,e\,(\frac{an}{c}) g (n) \\
  &= 2\pi i^k \frac{\eta}{c \sqrt r} \sum^\infty_1 \lambda_G (n)\,e\,
(\frac{- \overline{ar} n}{c}) \int^\infty_0 g(x) J_{k-1} \biggl(\frac{4\pi}{c} \sqrt
{\frac{nx}{r}} \biggr) dx \tag3.14
\endalign
$$ where $\lambda_G(n)$ are given by (2.28), (2.29) with $r =
q/(c,q), s = (c,q)$ and $$ \eta = \chi_s (a) \chi_r (-c) \eta_F
(r,s).  \tag3.15 $$ Here $\eta_F (r,s)$ depends only on $r,s,F$
and  $|\eta_F (r,s)| = 1$.
\endproclaim
\medskip
\demo{Proof}  The result follows by applying (3.11) for $A(z) = F(z)$ and
$B(z) = F_{|\omega}(z)$ with
$$
\omega = \left (\matrix a\sqrt r & b\, / \sqrt r\\ c\sqrt r & d\sqrt r
\endmatrix \right )  \tag 3.16
$$
where $b,d$ are integers such that $ad r - bc =1$.  Note that $d
\equiv {\overline{ar}}$(mod $c$) so the corresponding pseudo-eigenvalue
(2.31) is equal to (3.15).
\medskip
Next we apply (3.11) for the Eisenstein series $E_\fk(z)$ with
cusp $\fk \sim 1/v$ and the $\omega$ given by (3.16).  In this
case (2.38) holds with $\fk^\ast \sim 1/v^\ast = (r,v)^2/rv$,
where $r = q/(c,q), s = (c,q)$.  The corresponding
pseudo-eigenvalue (2.37) becomes $$ \eta = \chi_r
(\frac{c}{(c,v)}) \chi_s (\frac{av}{(c,v)}) \chi_{v/(c,v)} (-1).
\tag 3.17 $$ because $(s,v) = (c,v)$ and $(r,v) = v/(c,v)$.  We
have $q = vw = v^\ast w^\ast$ with $$ \aligned v^\ast &=
(c,v)w/(c,w)\\ w^\ast &= (c,w)v/(c,v).\endaligned  \tag 3.18 $$ We
 introduce the twisted divisor function $$ \tau(n;\chi_v,\chi_w) = \sum_{n_1n_2 = n} \chi_v(n_1)
\chi_w(n_2) \tag 3.19 $$ for any $n \geqslant 1$.   Therefore
$\tau(n;\chi_v,\chi_w)$ and $\tau(n;\chi_{v^\ast},\chi_{w^\ast})$
are the Hecke eigenvalues $\lambda_{\frak a}(n)$ and
$\lambda_{\fk^\ast}(n)$ for $E_\fk(z)$ and $E_{\fk^\ast}(z)$,
respectively.  These are proportional to the Fourier coefficients
of $E_\fk(z)$ and $E_{\fk^\ast}(z)$ with the factor $\bar \var_v
2i/h$ and $\bar \var_{v^\ast} 2i/h$, respectively (see (2.11)). By
(3.11) with $k=1$ and $c\sqrt r$ in place of $c$ we obtain $$
\align \sum^\infty_1 \tau(n;\chi_v,\chi_w)\,e\,(\frac{an}{c})g(n)
&= \var_v\bar \var_{v^\ast} \frac{2\pi i \eta}{c \sqrt r} \biggl\{
\frac{h}{2} \int^\infty_0 g(x) dx\\ &+ \sum^\infty_1
\tau(n;\chi_{v^\ast},\chi_{w^\ast})\,e\,
(\frac{-\overline{ar}n}{c}) \int^\infty_0 g(x) J_0 \biggl( \frac{4
\pi}{c} \sqrt {\frac{nx}{r}} \biggr) dx \biggr\}
\endalign
$$ where the leading term $\fc h \hat g (0)$ appears  only if
$\fk^\ast \sim \infty$, or ${\frak a}^* \sim 0$.  Note that $\bar
\var_{v^\ast} = \bar \var_{vr} = \bar \var_{-vs} = -i \var_{vs}$
so the factor $\sigma = i \eta \var_v \bar \var_{v^\ast}$ becomes
$$ \sigma = \var_v \var_{vs} \chi_s (\frac{av}{(c,v)}) \chi_r
(\frac{c}{(c,v)}) \chi_{(r,v)} (-1).  \tag 3.20 $$
\medskip
Before stating the final result we simplify the leading term.  We
have $\fk^\ast \sim \infty \Leftrightarrow (c,q) = v$ in which
case $r = w$ and $s = v$ so that $$ \sigma = \var_v \chi_v (a)
\chi_w (\frac{c}{v}) = \chi_v (a) \chi_w (\frac{c}{v}) \tau
(\chi_v)/\sqrt v. $$ Similarly $\fk^\ast \sim 0 \Leftrightarrow
(c,q) = w$ in which case $r = v$ and $s = w$ . Thus, $\sigma =
\var_v \var_q \chi_w(av) \chi_v(-c)$.  By the reciprocity law
$\chi_w(v) \chi_v(w) = 1$, because $vw = q \equiv -1$(mod $4$).
Moreover $\var_v \var_q \chi_v(-1) = i \bar \var_v = \var_{-v} =
\var_w$ so that $$ \sigma = \var_w \chi_w(a) \chi_v(\frac{c}{w}) =
\chi_w(a) \chi_v(\frac{c}{w}) \tau(\chi_w)/\sqrt w. $$ Hence we
conclude by the class number formula $\pi h = \sqrt q L(1,\chi)$
the following result.
\enddemo
\medskip
\proclaim{Proposition 3.3}  Let $q$ be squarefree, $q \equiv -1$(mod
$4$).
Let $q = vw$, and let $\chi_q = \chi_v \chi_w$ be the
corresponding real characters.  Let $c \geqslant 1$ and $(a,c) = 1$.  Let $q
= v^\ast w^\ast$ be given by (3.18), and let $\chi_q = \chi_{v^\ast}
\chi_{w^\ast}$ be the corresponding real characters.  Then for any
smooth function $g(x)$ compactly supported on $\Bbb R^+$ we have
$$
\align
\sum^\infty_1 \tau &(n;\chi_v,\chi_w)\,e\,(\frac{an}{c})g(n) \tag3.21 \\
  &= \{\chi_v(a)\chi_w(\frac{c}{v}) \tau(\chi_v) +
\chi_w(a)\chi_v(\frac{c}{w})
\tau(\chi_w)\} \frac{L(1,\chi)}{c} \int^\infty_0 g(x) dx\\
  &+ \frac{2\pi \sigma}{c\sqrt r} \sum^\infty_1 \tau
(n;\chi_{v^\ast},\chi_{w^\ast})\,e\, (\frac{-{\overline{ar}}n}{c})
\int^\infty_0 g(x) J_0 \biggl(\frac{4 \pi}{c} \sqrt {\frac{nx}{r}}
\biggr) dx
  \endalign
$$
where $\sigma$ is given by (3.20) with $r = q/(c,q)$ and $s =
(c,q)$.
\endproclaim
\medskip
\remark{Remark}.  In the leading term of (3.21) we use the popular
convention that $\chi(z)$ is zero if $z$ is not an integer.  Therefore
the leading term vanishes unless $v|c$ and $(c,w) = 1$, or $w|c$ and
$(c,v) = 1$.
\endremark
\medskip
\subhead 4.  Convolution Sums\endsubhead
\medskip
In this section we shall evaluate asymptotically sums of type $$
{\Cal B}(h) = \underset {m-n = h}\to {\dsize\sum\sum} \lambda(m)
\bar \lambda(n) g(m) \bar g(n) \tag4.1 $$ where $\lambda(n)$ are
eigenvalues of the Hecke operators $T_n$ for one of the
holomorphic automorphic forms of level $q$ and character $\chi_q$
which were considered in the last two sections.  Here $h \neq 0$
is a fixed integer and $g(x)$ is a cut-off function which is
smooth and compactly supported on $\Bbb R^+$.  Naturally one can
treat the convolution sum (4.1) by spectral methods using an
appropriate Poincar\'e series.  Indeed this method (the
Rankin-Selberg method) has been applied by various authors in many
cases of cusp forms, but we did not find satisfying results
\footnote{See the note added at the end of Section 1}.  Moreover
the spectral methods can also be applied to the Eisenstein series,
which case is important for us. However, the spectral methods are
rather complicated and one has to provide a lot of background
material.  Therefore in this section we use arguments from the
circle method of Kloosterman, because they yield the results
faster, more general and of great uniformity with respect to the
shift $h$.

For clarity we present the arguments in an axiomatic setting.
Literally speaking we do not assert that the coefficients
$\lambda(n)$ come from an automorphic form.  All we need is an
appropriate summation formula for $$ S(\alpha) = \sum^\infty_1
\lambda(n)\,e\,(\alpha n) g(n)  \tag 4.2 $$ at rational points. We
assume that for any $c \geqslant 1$ and $(a,c) = 1$ one has the
expansion $$ S(\frac {a}{c}) = \sum^\infty_{m=0} \psi_m
(a)\,e\,(\frac{\bar a}{c} l_m) \int^\infty_0 g(x) k_m(x) dx  \tag
4.3 $$ where $a\bar a \equiv 1$(mod $c$), $\psi_m (a)$ are
periodic in $a$, say of a fixed period $q, l_m$ are integers, and
$k_m(x)$ are smooth functions.  We do allow $\psi_m(a), l_m$ and
$k_m(x)$ to depend on $c$.  However the frequencies $l_m$ and the
kernels $k_m(x)$ cannot depend on $a$.  If $m \geqslant 0$ we
require the coefficients $\psi_m(a)$ to satisfy $$ |\psi_m(a)|
\leqslant A \tau(m) c^{-1}  \tag 4.4 $$ where $A$ is a constant,
$A \geqslant 1$.  Next we assume that the Fourier transform of
$g_m(x) = g(x) k_m(x)$ satisfies $$ |\hat g_m(\alpha)| \leqslant
Bc Cm^{-5/4} \quad \text{ if } |\alpha| \leqslant (cC)^{-1}  \tag
4.5 $$ for all $1 \leqslant c \leqslant C$, where $B \geqslant 1$
is a constant and $C \geqslant 2$ is a fixed number (a quite large
number which will be chosen optimally in applications). For $m =
0$ we need more precise conditions.  We assume that $$ l_0 =
0,\qquad k_0(x) = 1  \tag 4.6 $$ and the absolute value of
$\psi_0(a)$ does not depend on $a$, say $$ |\psi_0(a)| = p(c)
\leqslant A c^{-1}. \tag 4.7 $$ Finally we assume that $$ \align
\int^\infty_{-\infty} |\hat g(\alpha)| d\alpha \leqslant B \tag4.8
\\ \int^\infty_{-\infty} |\alpha||\hat g(\alpha)|^2 d\alpha
\leqslant B^2. \tag4.9
\endalign
$$
\medskip
Now we are ready to estimate the convolution sum (4.1).  We begin
by the following formula for the zero detector in $\Bbb Z$ $$ 2
\ssum \Sb c \leqslant C < d \leqslant c + C\\ (c,d) = 1 \endSb
\int^{1/cd}_0 \cos (2\pi n(\frac{a}{c} - \alpha))d\alpha =\cases 1
&\text{ if } n = 0\\ 0 &\text{ if } n \neq 0 \endcases  \tag4.10
$$ where $ad \equiv 1$ (mod $c$) and $C \geqslant 2$ is at our
disposal (see Proposition 11.1 of [I]).  Hence we get $$ {\Cal
B}(h) = \ssum \Sb c \leqslant C < d \leqslant  c + C\\ (c,d) = 1
\endSb \int^{1/cd}_{-1/cd} |S(\frac{a^\prime}{c} -
\alpha)|^2\,e\,(h(\alpha -\frac{a^\prime}{c}))d \alpha $$ where
$a'$(mod $c$) is determined by $a^\prime d \equiv \text{sign }
\alpha$ (mod $c$).  We rearrange this sum of integrals as follows
$$ {\Cal B}(h) = \dsize \sum_{c \leqslant C} \int^{1/cC}_{-1/cC} e
(\alpha h) V_c (\alpha) d \alpha  \tag 4.11 $$ where $$ V_c
(\alpha) = {\dsize \sum_{d \in I}}^* e(- \frac{a^\prime h}{c})
|S(\frac{a^\prime}{c} - \alpha)|^2 \tag4.12 $$ and $d$ runs over
integers prime to $c$ in the interval. $$ I = (C,\text{min}\{c +
C, \frac{1}{|\alpha|c}\}]. \tag 4.13 $$ Note that $I = (C, c + C]$
has length exactly $c$ if $|\alpha| \leqslant c^{-1}(c + C)^{-1}$,
and $I = (C,1/|\alpha|c]$ is shorter than $c$ if $c^{-1}(c +
C)^{-1} < |\alpha| \leqslant c^{-1}C^{-1}$.

Suppose $\alpha > 0$, the case $\alpha < 0$ is similar.  By the
summation formula (4.3) we have $$ S(\frac{a}{c} - \alpha) =
\sum^\infty_{m = 0} \psi_m (a)\,e\, (\frac{d}{c} l_m) \hat
g_m(\alpha). $$ Inserting this into (4.12) and changing the order
of summation we get $$ V_c(\alpha) = \sum_{m_1}\sum_{m_2} \hat
g_{m_1}(\alpha) {\bar{\hat g}_{m_2}}(\alpha) {\sum_{d \in I}}^*
\psi_{m_1} (a) \bar \psi_{m_2} (a)\,e\,(\frac{d}{c}(l_{m_1} -
l_{m_2})-\frac{ah}{c}). $$ Recall that $\psi_m(a)$ are periodic of
period $q$.  Splitting into residue classes $d \equiv \delta$ (mod
$q$) we obtain incomplete Kloostermann sums for which Weil's bound
yields $$ \sumstar \Sb d \in I ,\ d \equiv \delta (q) \endSb
\,e\,(\frac{dl - ah}{c}) \ll (h,c)^{\frac 1 2} c^{\frac 1 2}
\tau(c) \log C.  \tag 4.14 $$ We apply this result with $l =
l_{m_1} - l_{m_2}$ for all terms, except for $m_1 = m_2 = 0$ in
the range $|\alpha| < c^{-1}(c + C)^{-1}$, i.e. when the interval
(4.13) has length $c$.  We derive for any $\alpha$ $$ \align
V_c(\alpha) &= |\hat g(\alpha) p (c)|^2 \biggl\{ \sumstar_{C < d
\leqslant c + C}\, e\,(-\frac{ah}{c}) + O(|\alpha|c C(h,c)^{\frac
1 2} c^{\frac 1 2} \tau (c) \log C \biggr\} \\
  &+ O(A^2 q(h,c)^{\frac 1 2}c^{\frac 3 2}
\tau(c) (\log C) \biggl(|\hat g(\alpha)| + \sum^\infty_1
\tau(m)|\hat {g_m}(\alpha)|\biggr) \biggl(\sum^\infty_1
\tau(m)|\hat g_m(\alpha)|) \biggr) \endalign $$ by (4.4), (4.7)
and (4.14).  In the leading term we get exact Ramanujan sum $$
r_c(h) = \sumstar_{d(\text{mod }c)} \,e\,(\frac{dh}{c}). \tag4.15
$$ The same estimates hold for $\alpha < 0$.  Adding these results
we get by (4.11) $$ {\Cal B}(h) = \sum_{c \leqslant C}
r_c(h)p(c)^2 \int^{1/cC}_{-1/cC}\,e\,(\alpha h) |\hat g(\alpha)|^2
d\alpha + R $$ where $$ \align R &\ll A^2 \biggl( \int
|\alpha||\hat g(\alpha)|^2 d\alpha \biggr) \biggl( \sum_{c
\leqslant C} (h,c)^\frac 1 2 c^{-\frac 1 2} \tau(c)) C \log C \\
  &+ A^2 q \sum_{c \leqslant C} (h,c)^\frac 1 2 c^{-\fc 3 2}\tau(c)
\biggl( \int |\hat g(\alpha)|d\alpha + B \biggr) BcC \log C.
\endalign $$ To estimate $R$ we apply (4.8) and (4.9) getting $$
{\Cal B}(h) = \sum_{c \leqslant C} r_c(h)p(c)^2
\int^{1/cC}_{-1/cC}\,e\,(\alpha h)|\hat g(\alpha)|^2 d\alpha +
O(\tau(h)q A^2 B^2 C^{\frac 3 2}(\log C)^2). $$ In the leading
term we extend the integration to all $\alpha \in \Bbb R$ at the
cost of an error term which is already present.  Then we derive by
the Plancherel theorem that $$ \int^\infty_{-\infty}\,e\,(\alpha
h)|\hat g(\alpha)|^2 d \alpha = \int^\infty_0 g (x + h) \bar
g(x)dx. $$ Finally we extend the summation over $c \leqslant C$ to
all $c$ getting $$ \sum_{c \leqslant C} r_c(h)p(c)^2 = \sigma(h) +
O(\tau(h) A^2 C^{-1}) $$ where $\sigma(h)$ is the infinite series
$$ \sigma(h) = \sum^\infty_{c = 1} r_c(h)p(c)^2.  \tag 4.16 $$ We
have established the following
\medskip
\proclaim{Theorem 4.1}  Suppose the conditions (4.3) - (4.9) hold.
Then for any integer $h \neq 0$ the sum (4.1) satisfies
$$
\align
{\Cal B}(h) &= \{\sigma(h)
+ O(\tau(h) A^2 C^{-1})\} \int g (x + h)\bar g(x)dx \tag4.17 \\
  &+ O (\tau(h) q A^2 B^2 C^{\frac 3 2} (\log C)^2)  \endalign
$$
where $\sigma(h)$ is given by (4.16) and the implied constant is
absolute.
\endproclaim
\medskip
\remark{Remarks}  In the proof of Theorem 4.1 we assumed tacitly (just
to simplify notation) that the cut-off functions $g(x), \bar g(x)$ are
complex conjugate.  However the formula (4.17) holds (by obvious
alterations in the arguments) for any pair $g(x),\bar g(x)$, provided
both functions satisfy the same relevant conditions.  In forthcoming
applications we shall have two functions $g_1(x), g_2(x)$ supported in
[X,2X] with $X \geqslant \frac 1 2$ such that
$$
x^\nu|g^{(\nu)}_j (x)| \leqslant 1, \quad \text{ if } \nu = 0,1,2,  \tag 4.18
$$
for $j = 1,2$.  For such functions (4.8) and (4.9) hold with $B = 1$.
Moreover we shall be able to verify (4.5) with $C = 2 \sqrt {q X}$ and
some constant $B \geqslant 1$.  Therefore we state the following  \endremark
\medskip
\proclaim{Corollary 4.2}  Suppose the conditions (4.3) - (4.7) hold
for the arithmetic function $\lambda(n)$ and for the cut-off functions
$g_1(x), g_2(x)$ supported in [X,2X] with derivatives satisfying
(4.18).  Precisely let (4.5) hold with $C = 2\sqrt {q X}$.  Then for any
integer $h \neq 0$
$$
\align
{\Cal B}(h) &= \sum_{m - n = h} \lambda(m) \bar\lambda(n)g_1(m)g_2(n)
\tag 4.19 \\
  &= \sigma(h) \int g_1(x + h)g_2(x)dx + O(\tau(h)(q AB)^2
X^{\frac 3 4} (\log 3X)^2)
\endalign
$$
where $\sigma(h)$ is given by (4.16) and the implied constant is
absolute. \endproclaim
\medskip
\remark{Remarks}  It would be convenient for applications to have
a formula for ${\Cal B}(h)$ with the cut-off functions $g(x)$
smooth on $\Bbb R^+$ such that $$ x^\nu|g^{(\nu)}(x)| \leqslant (1
+ \frac x X)^{-4}, \quad \text{ if } \nu = 0,1,2,  \tag4.20 $$
rather than being supported in the dyadic segment [X,2X].
Unfortunately for such functions the condition (4.5) may not be
easily verifiable with a reasonable value of $C$ (the optimal $C$
should be of the order of $\sqrt X$).  Nevertheless we shall be
able to derive results for functions satisfying (4.20) by applying
a smooth partition of unity,  but not at the current position (see
how we justify (6.27)).
\endremark
\medskip
In principle our analysis (the Kloosterman circle method) works
also for $h = 0$, but, of course, giving somewhat different main
term.  In fact the resulting error term is better,  because the
estimate for the incomplete Kloosterman sum (4.14) is replaced by
a stronger bound for Ramanujan sum.  Rather than repeating and
modifying the former arguments we shall derive an asymptotic
formula for ${\Cal B}(0)$ directly using the Rankin - Selberg zeta
function,  see (6.43).
\medskip
Now we apply Corollary 4.2 for the $\lambda (n)$'s which are Hecke
eigenvalues of a holomorphic automorphic form of weight $k
\geqslant 1$, level $q$ and the real character $\chi_q$ of
conductor $q$.  As in Section 3 we assume that $q$ is odd, so $q$
is squarefree and $q \equiv 2k+1$ (mod $4$) by the consistency
condition $\chi_q (-1) = (-1)^k$.  This form is either a primitive
cusp form, or the holomorphic Eisenstein series $E_{\fk}(z) =
y^{-\frac 1 2} E_{\fk} (z, \frac 1 2)$ of weight $k = 1$ for a
cusp $\fk \sim 1/v$ with $vw = q$.  In the latter case the Hecke
eigenvalues are (see (3.19)) $$ \tau(n;\chi_v,\chi_w) = \dsize
\sum_{n_1n_2 = n} \chi_v(n_1) \chi_w(n_2). $$ The summation
formula (4.3) holds by Proposition 3.2 for cusp forms, or
Proposition 3.3 for the Eisenstein series.  In either case we have
$$ |\psi_m(a)| \leqslant 2\pi \tau(m) c^{-1}, \quad \text{ if } m
\geqslant 1. $$ and $|\psi_0(a)| = p(c)$ does not depend on $a$.
In fact $p(c) = 0$, except for the Eisenstein series $E_{\fk}(z)$
with $\fk \sim 1/v$, in which case we have $$ p(c) =
\frac{L(1,\chi)}{c} \cases \sqrt v &\text{ if } (c,q) = v\\ \sqrt
w &\text{ if } (c,q) = w \endcases  \tag4.21 $$ and $p(c) = 0$
otherwise.  In every case (cusp forms or Eisenstein series) the
summation formula holds with the kernel $$ k_m(x) = J_{k - 1}
\biggl( \frac{4\pi}{c} \sqrt {\frac{mx}{r}} \biggr) $$ where $r =
q/(c,q)$.  Note that for $k = 1$ we have $k_0(x) = 1$ as required
by (4.6).  Hence the Fourier transform of $g_m(x)$ is $$ \hat
g_m(\alpha) = \int g(x)\,e\,(\alpha x) J_{k - 1} \biggl(
\frac{4\pi}{c} \sqrt {\frac{mx}{r}} \biggr) dx.  \tag 4.22 $$ Note
that the Bessel function can be written as $$ J_{k - 1}(2\pi y) =
W(y)\,e\,(y) + \overline{W}(y)\,e\,(-y) $$ where $W(y)$ is a
smooth non-oscillatory function whose derivatives satisfy
 $$
y^\nu\, W^{(\nu)} (y) \ll k^2 y^{-\frac 1 2} \quad \text{ if } \nu
= 0,1,2. $$ Let g$(x)$ be a smooth function supported on [X,2X]
with $X \geqslant \frac 1 2$ such that $$ x^\nu |g^{(\nu)}(x)|
\leqslant 1, \quad \text{ if } \nu = 0,1,2. \tag 4.23 $$ We choose
$C = 2\sqrt {q X}$ so there is no stationary point in the Fourier
integral (4.22) if $1 \leqslant c \leqslant C$ and $|\alpha|c C
\leqslant 1$. Therefore integrating by parts two times we derive
$$ \hat g_m(\alpha) \ll k^2 X(c^2r/m X)^{\frac 5 4} \ll k^2
q^{\frac 3 2} cC m^{-\frac 5 4} $$ for $ |\alpha| \leqslant
(cC)^{-1}$, where the implied constant is absolute. Next we derive
from (4.23) by partial integration that $$ \hat g (\alpha) = \int
g(x)\,e\,(-\alpha x) dx \ll X (1 + |\alpha| X)^{-2}.  \tag 4.24 $$
Hence $$ \int |\hat g(\alpha| d\alpha \ll 1, \qquad
\int|\alpha||\hat g(\alpha)|^2 d\alpha \ll 1. $$ The above
estimates verify the conditions of Corollary 4.2 with $A = 2 \sqrt
q L(1,\chi), B = k^2 q^{\frac{3}{2}}$ and $C = 2 \sqrt {qX}$.
Hence we obtain the following two theorems.
\medskip
\proclaim{Theorem 4.3}  Let $\lambda_F(n)$ be the eigenvalues of a
primitive cusp form $F \in S_k(\Gamma_0(q), \chi_q)$ (recall that
$k \geqslant 1$ and $q$ is squarefree, $q \equiv 2k + 1$ (mod
$4$)).  Then for any integer $h \neq 0$ and for any smooth
functions $g_1(x), g_2(x)$ supported in [X,2X], $X \geqslant \frac
1 2$, with derivatives satisfying (4.23) we have $$ \underset {m -
n = h}\to {\dsize\sum\sum} \lambda_F(m) \bar \lambda_F(n) g_1(m)
g_2(n) \ll \tau(h)q^6 k^4 X^{\frac 3 4}(\log 3X)^2  \tag 4.25 $$
where the implied constant is absolute.
\endproclaim
\medskip
\proclaim{Theorem 4.4}  Let $q$ be squarefree, $q \equiv -1$(mod
$4$). Let $u w = q$ and $\tau(n; \chi_v,\chi_w)$ be the twisted
divisor function by the corresponding characters $\chi_v\chi_w =
\chi_q$ (see (3.19)). Then for any integer $h \neq 0$ and for any
smooth functions $g_1(x), g_2(x)$ supported in [X,2X], $X
\geqslant \frac 1 2$ with derivatives satisfying (4.23) we have $$
\align \underset {m - n = h}\to {\dsize\sum\sum}
&\tau(m;\chi_v,\chi_w) \tau(n;\chi_v,\chi_w) g_1(m) g_2(n) \tag
4.26 \\
  &= \sigma(h) \int g_1(x + h) g_2(x) dx + O(\tau(h) q^6 X^{\frac 3 4} (\log
3X)^2)
\endalign
$$
where $\sigma(h)$ is the infinite series (4.16) with $p(c)$ given by
(4.21), and the implied constant is absolute. \endproclaim
\medskip
\remark{Remark}  We emphasize that the estimates in the above
theorems are uniform in every parameter.  \endremark
\medskip
We conclude this section by computing $\sigma(h)$.  We have $$
\sigma(h) = \{ \sum_{(c,q) = v} \frac{v}{c^2} r_c(h) + \sum_{(c,q)
= w}\frac{w}{c^2} r_c(h) \} L^2(1,\chi)  \tag 4.27 $$ where
$r_c(h)$ is the Ramanujan sum.  Since $r_c(h)$ is multiplicative
in $c$ we get $$ \sigma(h) = \{\sum_{c|v^\infty}
\frac{r_{cv}(h)}{c^2v} + \sum_{c|w^\infty}
\frac{r_{cw}(h)}{c^2w}\}(\sum_{(c,q) = 1} \frac{r_c(h)}{c^2})
L^2(1,\chi) \tag 4.28 $$ where $$ \sum_{(c,q) = 1}
\frac{r_c(h)}{c^2} = \frac{\zeta_q(2)}{\zeta (2)} \dsize \sum \Sb
d|h \\ (d,q) = 1 \endSb d^{-1}.  \tag 4.29 $$ Moreover using the
formula $$ r_c(h) = \sum_{d|(c,h)} d \mu (c/d)  \tag 4.30 $$ one
can show that $$ \sum_{c|v^\infty} \frac{r_{cv}(h)}{c^2v} = \mu
(\frac{v}{(h,v)}) \frac{(h,v)}{v} \prod_{p|(h,v)} (1-
\frac{1}{p^\alpha} - \frac{1}{p^{\alpha + 1}}) \tag 4.31 $$ where
$p^\alpha || h$.  Gathering the above results we arrive at $$
\sigma(h) = \biggl\{ \mu(\frac{v}{(h,v)}) \frac{(h,v)}{v}
\underset {p|(h,v)}\to \prod (1 - \frac{1}{p^\alpha} -
\frac{1}{p^{\alpha + 1}}) + (v \to w) \biggr\}
\frac{\zeta_q(2)}{\zeta(2)} (\dsize \sum \Sb d|h \\ (d,q) = 1
\endSb \frac 1 d) L^2(1,\chi). \tag 4.32 $$ In applications we
shall appeal to the zeta-function of the $\sigma(h)$ $$ Z(s) =
\sum^\infty_1 \sigma(h) h^{-s}  \tag 4.33 $$ (note that $\sigma(h)
= \sigma(-h)$ because the Ramanujan sums are even in $h$).  Using
(4.27) and (4.30) one derives $$ \align Z(s) = \biggl\{ \frac 1 v
\underset {p|v}\to \prod (1- \frac{1}{p^{s-1}}) \underset {p|w}\to
\prod (1-\frac{1}{p^{s+1}}) + \frac 1 w \underset {p|w}\to \prod
(&1- \frac{1}{p^{s-1}}) \underset {p|v}\to \prod
(1-\frac{1}{p^{s+1}}) \biggr\}  \\
  & \frac{\zeta_q(2)}{\zeta(2)} \zeta(s) \zeta(s + 1)L^2(1,\chi).  \tag4.34
\endalign
$$

\medskip
Note that $Z(s)$ has no pole at $s = 1$ except for $v = 1$ or $w =
1$, i.e. if the cusp is at $\infty$ or $0$.  In these cases the
residue is $$ \underset{s=1} \to {res}\,\, Z(s) = L^2(1,\chi),
\quad \text{ if } v = 1 \text{ or } w = 1. \tag 4.35 $$ The only
other pole of $Z(s)$ is at $s = 0$ with residue $$ \underset{s =
0} \to {res}\,\, Z(s) = (\mu (v) + \mu (w)) \frac{q}{\nu (q)}
\frac{\zeta(0)}{\zeta(2)} L^2(1,\chi)  \tag 4.36 $$ where $\nu(q)$
is the multiplicative function with $\nu(p) = p +1$ (see (2.1)).
For curiosity we note that this residue vanishes if $\nu(q) = -1$,
for example if $q$ is prime.
\medskip
\subhead 5. Point to Integral Mean-Values of Dirichlet's
Series\endsubhead
\medskip
Our objective is to estimate a Dirichlet series $$ A(s) =
\sum^\infty_{n = 1} a_n n^{-s}  \tag 5.1 $$ on average with
respect to well-spaced points $s$.  In this section we transform
the problem to that for a corresponding integral in $s$. The
procedure is well-known and there is a variety of tools in the
literature, just to mention the original one by P.X. Gallagher
[G]. However the published results when applied directly to our
series do not always produce the desired effects.  What we need
are integrals which can be treated further by quite delicate
analysis in the off-diagonal range.  For this reason we cannot
afford to contaminate the coefficients $a_n$ by wild test
functions nor by sharp cuts. Therefore, rather than modifying the
existing results, we shall develop the desired transformations
from scratch.
\medskip
\proclaim{Lemma 5.1}  Let $a_n$ be any sequence of complex numbers
such that
$$
\sum_n |a_n| < \infty.  \tag 5.2
$$
Let $f(x)$ be a function of ${\Cal C}^1$ class on $[1,\infty)$ such that
$$
c_f = \int^\infty_1 (x^{-1} |f(x)|^2 + x |f^\prime(x)|^2) dx \,\,<
\infty.  \tag 5.3
$$
Then we have
$$
|\sum_n a_n f(n)|^2\, \leqslant \frac{c_f}{\pi} \int^\infty_{-\infty}
|A(it)|^2 \frac{dt}{t^2 + 1}.  \tag 5.4
$$
\endproclaim
\medskip
\demo{Proof}  First we extend $f(x)$ to the segment $[0,1)$ by
setting $f(x) = x f(1)$.  Then we write $$ f(x) = \frac{1}{2\pi}
\int^\infty_{-\infty} h(t) x^{-it} dt $$ where $$ h(t) =
\int^\infty_0 f(x) x^{it - 1} dx $$ by Mellin (or Fourier)
inversion.  This gives us $$ \sum_n a_n f(n) = \frac{1}{2\pi}
\int^\infty_{-\infty} h(t) A(it) dt. $$ Hence by Cauchy-Schwarz
inequality $$ |\sum_n a_n f(n)|^2 \leqslant \frac{1}{4\pi^2}
\biggl( \int^\infty_{-\infty} |h(t)|^2 (t^2 + 1) dt \biggr)
\int^\infty_{-\infty} |A(it)|^2 \frac{dt}{t^2 + 1}. $$ By
Plancherel's theorem $$ \int^\infty_{-\infty} |h(t)|^2 dt = 2\pi
\int^\infty_0 x^{-1} |f(x)|^2 dx $$ $$ \int^\infty_{-\infty}
|h(t)|^2 t^2 dt = 2\pi \int^\infty_0 x |f^\prime(x)|^2 dx. $$
Hence we obtain (5.4) with the constant $$ c^\ast_f = \frac{1}{2}
\int^\infty_0 (x^{-1} |f(x)|^2 + x|f^\prime(x)|^2) dx $$ in place
of $c_f$.  Here the integral over the segment $[0,1)$ equals
$|f(1)|^2$.  Moreover we have $$ \align f(1)^2 &= - \int^\infty_1
(f^2(x))^\prime dx = - 2 \int^\infty_1 f(x) f^\prime(x) dx\\
  &\leqslant \int^\infty_1 (x^{-1} |f(x)|^2 + x |f^\prime(x)|^2) dx. \endalign
$$ Hence $c^\ast_f \leqslant  c_f$ proving (5.4).
\enddemo
\medskip
\proclaim{Corollary 5.2}  Let the conditions be as in Lemma 5.1.  Then
for $\rho = \beta + i \gamma$ with $0 \leqslant \beta \leqslant
\frac{1}{2}$ we have
$$
|\sum_n a_n n^{-\rho}f(n)|^2 \,\leqslant \frac{2c_f}{\pi}
\int^\infty_{-\infty} |A(it)|^2 \frac{dt}{(t-\gamma)^2 + 1}.  \tag 5.5
$$
\endproclaim
\medskip
\demo{Proof}  Apply (5.4) for $a_n n^{- i\gamma}$ and $n^{-\beta}f(n)$ in
place of $a_n$ and $f(n)$.
\medskip
Let $\Cal R$ be a set of points $\rho_r = \beta_r + i\gamma_r$ for
$r = 1,2,...,R$ such that $$ 0 \leqslant \beta_r \leqslant
\frac{1}{2}  \tag 5.6 $$ $$ T \leqslant \gamma_r \leqslant 2T \tag
5.7 $$ $$ |\gamma_r - \gamma_{r^\prime}| \geqslant \delta \quad
\text{ if } r \neq r^\prime.  \tag 5.8  $$ Here $\delta, T$ are
fixed numbers with $0 < \delta \leqslant 1$ and $T \geqslant 2$.
Note that $R \leqslant 1 + \frac{T}{\delta} \leqslant
\frac{3T}{2\delta}$.
\medskip
Suppose for every $r = 1,2,...,R$ we have a function $f_r(x)$ of
${\Cal C}^1$ class on $[1,\infty)$ such that corresponding
integrals (5.3) are bounded.  Put $$ c = \underset r \to \max
\,\,c_{f_r}  \tag 5.9 $$  $$ G(t) = \sum_r ((t-\gamma_r)^2 +
1)^{-1}. \tag 5.10 $$ From (5.5) we get immediately $$
\sum_r|\sum_n a_n n^{-\rho_r}f_r(n)|^2 \leqslant \frac{2c}{\pi}
\int^\infty_{-\infty} |A(it)|^2 G(t) dt.  \tag 5.11 $$
\medskip
Now we are going to estimate $G(t)$. By the spacing condition (5.8)
we derive that
$$
\align
G(t) &\leqslant 1 + \sum^\infty_1 (\delta^2 n^2 + 1)^{-1}\\
  &\leqslant 1+ \int^\infty_0 (\delta^2 t^2 + 1)^{-1} dt \,= 1 +
\frac{\pi}{2\delta} < \frac{\pi}{\delta}.
\endalign
$$ If $t$ is far beyond the segment (5.7) we can do better. Indeed
if $t < \frac{T}{2}$ then $(t - \gamma_r)^2 \geqslant (t - T)^2
\geqslant \frac{1}{5}(t^2 + T^2)$ and if $t > 3T$ then $(t -
\gamma_r)^2 \geqslant (t -2T)^2 \geqslant \frac{1}{10}(t^2 +
T^2)$.  Hence in these ranges $$ G(t) \leqslant \frac{10R}{t^2 +
T^2}\,\leqslant\,\frac{15T}{\delta (t^2 + T^2)}. $$ Inserting
these estimates into (5.11) we get $$ \align \sum_r|\sum_n a_n
n^{-\rho_r}f_r(n)|^2 &\leqslant \frac{2c}{\delta} \int^{3T}_{T/2}
|A(it)|^2 dt  \tag 5.12  \\
  &+ \frac{10c}{\delta T} \int^\infty_{-\infty} |A(it)|^2 (1 +
\frac{t^2}{T^2})^{-1} dt.
\endalign
$$ The first integral is no larger than ten times of the second
one, so we have  $$ \sum_r|\sum_n a_n n^{-\rho_r} f_r(n)|^2
\leqslant 30{\frac{c}{\delta}} \int^\infty_{-\infty} |A(it)|^2 (1
+ \frac{t^2}{T^2})^{-1} dt \tag 5.13 $$ (later we shall do better
with the first integral). The last integral is exactly equal to $$
\int^\infty_{-\infty} |A(it)|^2 (1 + \frac{t^2}{T^2})^{-1} dt =
\pi T \sum_m \sum_n a_m \bar a_n \text{
min}\biggl(\frac{m}{n},\frac{n}{m}\biggr)^T. $$ Now assuming that
$$ G_1 = \sum^\infty_1 n |a_n|^2 \,\, < \infty \tag 5.14 $$ we
estimate as follows  $$ \sum_m \sum_n a_m \bar a_n \text{
min}\biggl(\frac{m}{n},\frac{n}{m}\biggr)^T\,\leqslant \sum_n
|a_n|^2 \sum_{m \geqslant n} (\frac{n}{m})^T, $$  $$ \sum_{m
\geqslant n} (\frac{n}{m})^T \leqslant 1 + \int^\infty_n
(\frac{n}{x})^T dx = 1 + \frac{n}{T-1} \leqslant 2(1 +
\frac{n}{T}). $$ Hence $$ \int^\infty_{-\infty} |A(it)|^2
(1+\frac{t^2}{T^2})^{-1} dt \leqslant 2\pi (TG + G_1)  \tag 5.15
$$ where $$ G = \sum_n |a_n|^2.  \tag 5.16 $$ Inserting (5.15)
into (5.13) we get
\enddemo
\medskip
\proclaim{Lemma 5.3}  Let $\rho_r$ and $f_r(x)$ be as above.  Suppose
the complex numbers $a_n$ satisty (5.2) and (5.14).  Then
$$
\sum_r |\sum_n a_n n^{-\rho_r}f_r(n)|^2 \leqslant 189 \frac{c}{\delta}(TG + G_1)
\tag 5.17
$$
where $G, G_1$ are defined by (5.16) and (5.14). \endproclaim
\medskip
The estimate (5.17) (nevermind the constant 189) is not
sufficiently strong when the range of coefficients $a_n$ exceeds
$T$. Having this case in mind we retain the first integral in
(5.12) and apply (5.15) only to the second one. Actually we
enlarge the first integral slightly while smoothing the
integration. Precisely we set $$ {\Cal A}(T) =
\int^\infty_{-\infty} K (\frac{t}{T}) |A(it)|^2 dt  \tag 5.18 $$
where $K(u)$ is a non-negative function on $\Bbb R$ such that
$K(u) \geqslant 1$ for $\frac{1}{2} \leqslant u \leqslant 3$.  We
obtain
\medskip
\proclaim{Proposition 5.4}  Let $\rho_r,f_r(x)$ and $a_n$ be as
above. Then $$ \sum_r|\sum_n n^{-\rho_r}f_r(n)|^2 \leqslant
\frac{2c}{\delta} {\Cal A}(T) + \frac{63c}{\delta T} (TG + G_1).
\tag 5.19 $$
\endproclaim
\medskip

\noindent
\subhead {6. Evaluation of ${\Cal A}(T)$} \endsubhead
\medskip
By (5.15) one gets the bound  $\Cal A(T) \ll TG + G_1$ which is
essentially best possible in general. In this section we evaluate
$\Cal A(T)$ more precisely  for special sequences $\Cal A =
(a_n)$.  We assume that the cut-off function $K(u)$ in the
integral (5.18) is continuous and symetric on $\Bbb R$ with $$
K(0) = 0. \tag6.1 $$ Moreover we assume that the cosine-Fourier
transform $$ L(v) = 2 \int^\infty_0 K(u) \text{cos}(uv)du \tag6.2
$$ has fast decaying derivatives, specifically $$ |L^{(j)}(v)|
\leqslant (1 + |v|)^{-4}, \qquad 0 \leqslant j \leqslant 5.
\tag6.3 $$ Clearly any smooth, symmetric and compactly supported
function on $\Bbb R \backslash \{0\}$ does satisfy the above
conditions up to a constant factor. We get $$ {\Cal A}(T) = T
\dsize \sum_m \, \dsize \sum_n \, a_m \bar a_n L(T \log
\frac{m}{n}). $$ Here $L(T \log\frac{m}{n})$ localizes the terms
close to the diagonal. Therefore, we arrange this double sum
according to the difference $m-n=h$ with the intention to treat
every partial sum $$ S(h) = \dsize \sum_{m-n=h} \, a_m \bar a_n
L(T \log\frac{m}{n}) \tag6.4 $$ separately.  Note that $S(-h) =
\overline{S(h)}$ so we have $$ {\Cal A}(T) = L(0)TG + 2T \text{ Re
} \dsize \sum_{h > 0} S(h). \tag6.5 $$ Recall that $G$ is given by
(5.16). Here the zero term comes from the diagonal $m=n$ ; we have
$S(0) = L(0)G$ and $$ L(0) = 2 \int^\infty_0 K(u)du. \tag6.6 $$
Let $h > 0$. Thinking of $h$ as being relatively small we use the
approximation $$ \log \frac{m}{n} = \log(1 + \frac{h}{n}) =
\frac{h}{n} + O(\frac{h^2}{n^2}) $$ to modify $L(T \log
\frac{m}{n})$ as follows $$ L(T \log \frac{m}{n}) =
L(\frac{hT}{n}) + O(\frac{1}{T} (1 + \frac{hT}{n})^{-2}). $$ The
contribution of the error term to $S(h)$, say $S'(h)$, satisfies
$$ \align S'(h) &\ll \frac{1}{T} \dsize \sum_{m-n=h} |a_ma_n|(1 +
\frac{hT}{n})^{-2} \\
  &\ll \frac{1}{T} \dsize \sum_{m-n=h} (|a_m|^2 + |a_n|^2)
(1 + \frac{hT}{n})^{-2}.
\endalign
$$ Hence the contribution of the error terms to ${\Cal A}(T)$, say
${\Cal A}'(T)$, satisfies $$ {\Cal A}'(T) \ll \underset{m > n}\to
{\sum \sum} (|a_m|^2 + |a_n|^2) (1 + \frac{m-n}{n} T)^{-2}. $$
Hence it follows that $$ {\Cal A}'(T) \ll G + T^{-1} G_1. \tag6.7
$$ We are left with $$ {\Cal A}(T) = L(0)TG + 2T \text{ Re }
\dsize \sum_{h > 0} S^*(h) + O(G + T^{-1}G_1) \tag6.8 $$ where
$S^*(h)$ is the modified sum $$ S^*(h) = \dsize \sum_n a_{n+h}
\,\,\bar a_n L(\frac{hT}{n}). \tag6.9 $$ We may estimate $S^*(h)$
trivially as follows $$ S^*(h) \ll \dsize \sum_n (|a_{n+h}|^2 +
|a_n|^2)(\frac{n}{hT})^2 \leqslant 2(hT)^{-2} G_2 \tag6.10 $$
subject to the condition $$ G_2 = \dsize \sum_n n^2 |a_n|^2 \,\, <
\infty. \tag6.11 $$ This estimate is quite useful for large $h$,
say $h \geqslant H$, where $H$ will be defined later. Inserting
(6.10) into (6.8) we get $$ {\Cal A}(T) = L(0)TG + 2T \text{ Re }
\dsize \sum_{0 < h \leqslant H} S^*(h) + O(G+T^{-1}G_1 + T^{-1}
H^{-1}G_2). \tag6.12 $$ Observe that the terms of (6.9) with $n
\leqslant 2h$ contribute less than $$ \dsize \sum_{n \leqslant 2h}
|a_{n+h}\,\, a_n L(\frac{hT}{n})| \ll (hT)^{-2} \dsize \sum_{n
\leqslant 3h} n^2 |a_n|^2. $$ Summing over $h$ we find that these
small terms contribute to ${\Cal A}(T)$ less than $$  T^{-1}
\dsize \sum_n n^2|a_n|^2 \dsize \sum_{3h \geqslant n} h^{-2} \ll
T^{-1} G_1 $$ which is absorbed by the error term already present
in (6.12).

Now we require that the coefficients $a_n$ are given by $$ a_n =
\lambda(n)a(n) \tag6.13 $$ where $\lambda(n)$ is a nice arithmetic
function and $a(y)$ is a smooth cut-off function.  We do not
restrict $a(y)$ to a dyadic segment, but for practical needs we
require only that $a(y)$ is a ${\Cal C}^2$ class function on $\Bbb
R^+$ such that $$ y^\nu|a^{(\nu)}(y)| \leqslant (1 +
\frac{y}{Y})^{-4}, \quad \text{ if } \quad \nu = 0,1,2, \tag6.14
$$ where $Y \geqslant 2$.  Concerning $\lambda(n)$ we assume that
it is bounded by the divisor function $$ |\lambda(n)| \leqslant
\tau(n). \tag6.15 $$ Therefore our coefficients are almost
bounded, precisely  $$ |a_n| \leqslant \tau(n)(1 +
\frac{n}{Y})^{-4}. \tag6.16 $$ Hence the series of $|a_n|^2$,
$n|a_n|^2$, $n^2|a_n|^2$ converge and satisfy $$ G \ll Y (\log
Y)^3, \quad G_1 \ll Y^2 (\log Y)^3, \quad G_2 \ll Y^3 ( \log Y)^3.
\tag6.17 $$

Moreover about $\lambda(n)$ we postulate that for every two smooth
functions $g_1(x),g_2(x)$ supported in $[X,2X]$ with $X \geqslant \frac12$
such that
$$
x^\nu|g^{(\nu)}_j(x)| \leqslant 1, \quad \text{ if } \quad \nu = 0,1,2 \tag6.18
$$
and for any $h \geqslant 1$ we have
$$
\align
\dsize \sum_{m-n=h} \lambda(m) \bar \lambda(n)g_1(m)g_2(n) &=
\sigma(h) \int g_1(x+h)g_2(x)dx \tag6.19 \\
  &+ O(B \tau(h) X^{\frac{3}{4}} (\log\, 3X)^2).
\endalign
$$ Here $\sigma(h)$ is another nice arithemetic function depending
on $\lambda$, $B$ is a positive constant depending on $\lambda$,
and the implied constant in the error term is absolute.  We assume
that $$ |\sigma(h)| \leqslant C \dsize \sum_{d|h} d^{-1}. \tag6.20
$$ In other words the generating series $$ Z(s) = \dsize
\sum_{h=1}^\infty \sigma(h)h^{-s} \tag6.21 $$ is majorized by $C
\zeta(s)\zeta(s+1)$.  More precisely we assume that $$ Z(s) =
\zeta(s)z(s) \tag6.22 $$ where $z(s)$ is holomorphic in Re $s
\geqslant 0$, except for a simple pole at $s=0$.  Suppose $$ |z(s)
- \frac{Z}{s}| \leqslant C(|s|+1) \tag6.23 $$ in the strip $0
\leqslant Re \, s \leqslant \frac{3}{2}$.  We do not exclude the
residue $Z=0$, and we assume $$ |Z| \leqslant C. \tag6.24 $$

For our primary example we let the  $\lambda(n) = \lambda_F(n)$ be
the Hecke eigenvalues of a primitive cusp form $F \in
S_k(\Gamma_0(q),\chi_q)$. In this case (6.15) is proved by P.
Deligne (the Ramanujan conjecture) and the formula (6.19) is
established in our Theorem 4.3 with $\sigma(h)=0$ and $B =
k^4q^6$. Therefore $Z=0$ and $C=0$.

Our second example is the Hecke eigenvalue $\lambda(n) =
\tau(n;\chi_v,\chi_w)$ of a holomorphic Eisenstein series of
weight $k=1$ and level $q$.  In this case (6.15) is obvious by
(3.19) and the formula (6.19) is established in our Theorem 4.4
with $\sigma(h)$ given by (4.32) and $B = q^6$. The generating
series $Z(s)$ is computed in (4.34) and the residue of $z(s)$ at
$s=0$ is $$ Z = - \frac{3}{\pi^2} (\mu(v) + \mu(w))
\frac{q}{\nu(q)} L^2(1,\chi_q) \tag6.25 $$ (see (4.36)).  In this
case the estimates (6.20), (6.23), (6.24) hold with $$ C \ll
\frac{\nu(q)}{q} L^2(1,\chi_q) \log\, q. \tag6.26 $$

Now we are ready to evaluate $S^*(h)$.  By (6.19) we derive $$
\align S^*(h) &= \sigma(h) \int a(y+h) \bar a(y) L(\frac{hT}{y})dy
 \tag6.27 \\
  &+ O(B\tau(h)Y^{\frac{3}{4}} (\log Y)^4) + T^{-2}h|\sigma(h)| + T^{-2}
h(\log\, 3h)^3).
\endalign
$$ Well, not immediately because $a(y)$ is not supported in a
dyadic segment.  However, using a smooth partition of unity with
constituents in $m,n$ supported in segments of type $[Y_1,\sqrt 2
Y_1],[Y_2,\sqrt 2 Y_2]$ respectively one can justify the
applicability of (6.19) as follows.  Indeed there is no question
when the two segments $[Y_1,\sqrt 2 Y_1],[Y_2,\sqrt 2 Y_2]$ are
equal or adjacent.  If these segments are separated then they
produce nothing from the sum nor from the integral in (6.27)
unless $Y_1,Y_2 \leqslant \sqrt 2 h$.  In this case we estimate
trivially by $$ \dsize \sum_{n \leqslant 2h} |a_{n+h}\,
a_n|(\frac{n}{hT})^2 \leqslant \dsize \sum_{n \leqslant 3h}
|a_n|^2(\frac{n}{hT})^2 \ll T^{-2} h(\log \,3h)^3 $$ which yields
the third error term in (6.27).  Moreover, the integral over $x
\leqslant 2h$ is estimated similarly by $$ \int^{2h}_0
|a(x+h)a(x)|(\frac{x}{hT})^2 dx \leqslant \int^{3h}_0 |a(y)|^2
(\frac{y}{hT})^2 dy \ll T^{-2} h $$ which yields the second error
term in (6.27).

Next we replace $a(y+h)$ in (6.27) by $a(y)$ with the difference
$O(Y^2 T^{-2} h^{-1} |\sigma(h)|)$.  We obtain
$$
\align
S^*(h) &= \sigma(h) \int |a(y)|^2 L(\frac{hT}{y})dy \tag6.28 \\
  &+ O(B \tau(h)(Y^{\frac{3}{4}} + Y^2 T^{-2} h^{-1} + T^{-2} h)(\log \,
hY)^4)
\endalign
$$
where the implied constant is absolute.  This is true for all $h \geqslant
1$, but we only use this for $1 \leqslant h \leqslant H$,
where $H$ will be chosen
later.  Introducing (6.28) into (6.12) we derive
$$
\align
{\Cal A}(T) &= L(0)TG + 2T \int |a(y)|^2 \biggl(\dsize \sum^H_{h=1}
\sigma(h) L(\frac{hT}{y})\biggr)dy \\
  &+ O(B(THY^{\frac{3}{4}} + T^{-1} Y^2 + T^{-1}H^2)(\log HY)^5) \\
  &+ O((Y + T^{-1} Y^2 + T^{-1} H^{-1} Y^3)(\log Y)^3).
\endalign
$$ Note that we can extend the sum over $1 \leqslant h \leqslant
H$ to the infinite series $$ D(v) = \dsize \sum^\infty_{h=1}
\sigma(h) L(hv) \tag6.29 $$ with $v = Ty^{-1}$, up to the error
term $O(T^{-1} H^{-1} Y^3 \log H)$ which is already present (the
last one).  Having done this we choose $$ H = B^{-\frac{1}{2}}
T^{-1} Y^{\frac{9}{8}} (\log Y)^{-1} \tag6.30 $$ (this choice
equalizes the first and the last error terms) getting $$ \align
{\Cal A}(T) &= L(0)TG + 2T \int |a(y)|^2 D(T/y)dy \tag6.31 \\
  &+ O(B^{\frac{1}{2}} Y^{\frac{15}{8}} \log^4 Y + T^{-3}
Y^{\frac{9}{4}} \log^3 Y).
\endalign
$$

Next we evaluate the series $D(v)$.  Let $M(s)$ be the Mellin
transform of $L(v)$ $$ M(s) = \int^\infty_0 L(v)v^{s-1} dv. $$
Integrating by parts we derive by (6.2) $$ sM(s) \ll (|s|+1)^{-4}.
\tag6.32 $$ Note that $$ M(1) = \int^\infty_0 L(v)dv = 2\pi K(0) =
0 $$ by our assumption (6.1).  Therefore the product $M(s)Z(s)$ is
holomorphic in the strip $0 < \sigma \leqslant \frac{3}{2}$ (no
pole at $s=1$) and $$ M(s)Z(s) \ll C|s|^{-2} $$ by (6.22), (6.23),
(6.24) and (6.26).  However for $s$ near zero we need a more
precise expansion.  To this end we use $$ \zeta(s) = \zeta(0) +
O(|s|), \qquad z(s) = \frac{Z}{s} + O(C), $$ and we derive an
expansion for $M(s)$ as follows $$ \align M(s) &= \frac{L(0)}{s} +
\int^1_0 (L(v) - L(0)) v^{s-1}dv + \int^\infty_1 L(v) v^{s-1} dv
\\
  &= \frac{L(0)}{s} + O(1).
\endalign
$$
From these expansions we get
$$
M(s)Z(s) = \frac{a}{s^2} + \frac{b}{s} + O(C)
$$
where $a = \zeta(0)L(0)Z = -\frac{1}{2} L(0)Z \ll C$ and $b \ll C$.
Combining both estimates we get
$$
M(s)Z(s) = \frac{a}{s^2} + \frac{b}{s(s+1)} + O(\frac{C}{|s|^2 +1}) \tag6.33
$$
uniformly in $0 < \sigma \leqslant \frac{3}{2}$.  Now we are ready to
evaluate $D(v)$.  We have
$$
\align
D(v) &= \frac{1}{2 \pi i} \int_{(\sigma)} M(s)Z(s)v^{-s} ds \\
  &= a \,\log^+ \frac{1}{v} + b \max(0,1-v) + O(Cv^{- \sigma}).
\endalign
$$
Hence we write
$$
D(v) = -\frac{1}{2} L(0)Z \log^+ \frac{1}{v} + D_0(v) \tag6.34
$$
with $D_0(v)$ a bounded function, specifically
$$
D_0(v) \ll C \tag6.35
$$
by letting $\sigma \rightarrow 0$ (recall the uniformity in $\sigma$).
Inserting (6.34) into (6.31) we conclude the above considerations by
the following
\bigskip

\proclaim{Theorem 6.1}  Let $K(u)$ be a continuous and symmetric
function on $\Bbb R$ with $K(0) = 0$ such that (6.2) holds.  Let
$\lambda(n)$ be an arithmetic function with $|\lambda(n)| \leqslant \tau(n)$
which satisfies the formula (6.19) with the surrounding conditions
(6.18) - (6.24).  Let $a(y)$ be a ${\Cal C}^2$ class function
on $\Bbb R^+$ such that (6.14) holds.  Then we have
$$
\align
\int^\infty_{- \infty} &K(\frac{t}{T})| \dsize \sum_n
a(n)\lambda(n)n^{-it}|^2 dt \tag6.36  \\
  &= \hat K(0)T \{G-Z \int^\infty_T |a(y)|^2 (\log \frac{y}{T})dy\} \\
  &+ 2T \int^\infty_0 |a(y)|^2 D_0(\frac{T}{y})dy
  + O(B^{\frac{1}{2}} Y^{\frac{15}{8}} (1 + T^{-3}
Y^{\frac{3}{8}})(\log Y)^4)
\endalign
$$
where
$$
G = \dsize \sum_n |a(n)\lambda(n)|^2 \tag6.37
$$
\mn
and $B,C,Z$ are the constants depending on $\lambda(n)$ given by the
postulated properties (6.18) - (6.24).  Moreover $D_0(v)$ is defined
by (6.34), so $D_0(v) \ll C$, the implied constants being absolute.
\endproclaim
\medskip

>From (6.36) one can derive a mean-value theorem for $A(s)$ on the line
Re $s = \frac{1}{2}$, however not without some loss in the error term.
We do it for $A(\frac{1}{2} + it)$ localized between $T$ and $Y$.
Precisely we get
\proclaim{Corollary 6.2}  Let the conditions be as in Theorem 6.1,
except for (6.14) which is now replaced by
$$
y^\nu |a^{(\nu)}(y)| \leqslant \biggl(
1 + \frac{y}{Y} + \frac{T}{y} \biggr)^{-4}, \quad
 \text{ if } \quad \nu = 0,1,2, \tag6.38
$$
where $T \leqslant Y \leqslant T^8$.  Then
$$
\align
\int^\infty_{-\infty} &K(\frac{t}{T})| \dsize \sum_n
a(n)\lambda(n)n^{\frac{1}{2}-it}|^2 dt \tag6.39  \\
  &= \hat K(0)T \{G-Z \int^\infty_T |a(y)|^2 (\log \frac{y}{T})
\frac{dy}{y}\} \\
  &+ 2T \int^\infty_0 |a(y)|^2 D_0 \bigl(\frac{T}{y}\bigr) \frac{dy}{y}
  + O(B^{\frac{1}{2}} T^{-1} Y^{\frac{15}{8}} \log^4 Y)
\endalign
$$
where $Z,D_0(v)$ and $B$ are as before, but
$$
G = \dsize \sum_n |a(n)\lambda(n)|^2 n^{-1}. \tag6.40
$$
\endproclaim
\bigskip

\demo{Proof}  Apply (6.36) for the function $a(y) \sqrt{T/y}$ in place
of $a(y)$.  This modified function satisfies (6.14) apart of an
absolute constant factor by virtue of (6.38).  Then divide the
resulting formula throughout by $T$.

Estimating all but the first term on the right side of (6.39) we
obtain
$$
\align
\int^\infty_{-\infty} K(\frac{t}{T})&|A(\tfrac{1}{2} + it)|^2 dt = \hat K(0)TG  \tag6.41 \\
  &+ O\bigl( T(|Z|\log \frac{Y}{T} +C)\log \frac{Y}{T} + B^{\frac{1}{2}}
T^{\frac{15}{16}} \log^4 T \bigr)
\endalign
$$
if $2T \leqslant Y \leqslant T^{\frac{31}{30}}$.  Moreover $G \ll (\log Y)^3
\log(Y/T)$ by (6.15).  But we are looking for a better estimate of $G$;
besides reducing by logarithms we want to see the implied constant.

We are most interested in Hecke eigenvalues $\lambda(n)$ of
automorphic forms associated with the imaginary quadratic field $K
= \Bbb Q(\sqrt{-q})$.  For these the Ramanujan bound (6.15) can be
improved to $$ |\lambda(n)| \leqslant \tau(n,\chi). \tag6.42 $$
The extremal case $\lambda(n) = \tau(n,\chi)$ comes from the
Eisenstein series $E_{\frak a}(z) = y^{-\frac{1}{2}} E_{\frak a}
(z,\fc)$ for cusps ${\frak a} =0$ or ${\frak a} = \infty$. The
zeta function of $\tau(n,\chi)$ is the $L$-function of $K = \Bbb
Q(\sqrt{-q})$, $$ L_K(s) = \dsize \sum_{0 \ne {\frak a} \subset
{\Cal O}_K} (N{\frak a})^{-s}  = \dsize \sum^\infty_{n=1}
\tau(n,\chi)n^{-s} = \zeta(s)L(s,\chi). $$ The Rankin-Selberg
$L$-function is $$ R_K(s) = \dsize \sum^\infty_{n=1}
\tau^2(n,\chi)n^{-s} = \zeta^2(s)L^2(s,\chi)\, \zeta(2s)^{-1}
\dsize \prod_{p|q} \bigl( 1 + \frac{1}{p^s} \bigr)^{-1}. \tag6.43
$$ This has the Taylor expansion $$ R_K(s) =
\frac{\alpha}{(s-1)^2} + \frac{\beta}{s-1} + \gamma + \ldots $$
with the polar coefficients given by $$ \alpha =
\frac{q}{\nu(q)}\, \frac{L^2(1,\chi)}{\zeta(2)} \tag6.44 $$ $$
\beta = \frac{q}{\nu(q)} \, \frac{L^2(1,\chi)}{\zeta(2)}\biggl[
2\frac{L'}{L}(1,\chi) + \gamma_1 + \dsize \sum_{p|q} \frac{\log
p}{p+1} \biggr]. \tag6.45 $$ Moreover $$ R_K(s) \ll
q^{\frac{1}{2}} |s|^{\frac{5}{6}}, \text{ if } \text{ Re } s
\geqslant \tfrac{1}{2}. $$ Let $g(y)$ be a function of ${\Cal
C}^2$ class on $\Bbb R^+$ such that $$ y^\nu|g^{(\nu)}(y)|
\leqslant \biggl(1 + \frac{y}{Y} + \frac{X}{y} \biggr)^{-1},
\text{ if } \nu = 0,1,2 \tag6.46 $$ where $Y \geqslant X > 0$. Let
$\check g(s)$ be the Mellin transform of $g(y)$, $$ \check g(s) =
\int^\infty_0 g(y)y^{s-1} dy = \int^\infty_0 g(y)(1+s\, \log\, y +
\ldots) \frac{dy}{y}. $$ By partial integration we get $$ \check
g(s) \ll |s|^{-2} (X^\sigma + Y^\sigma) \quad \text{ if } \quad
\sigma = \text{ Re } s = \pm \tfrac{1}{2}. $$ By contour
integration the sum $$ G = \dsize \sum_n \tau^2(n,\chi)
\frac{g(n)}{n} \tag6.47 $$ \mn is equal to $$ \align G &=
\frac{1}{2 \pi i} \int_{(\frac{1}{2})} \check g(s) R_K(s+1) ds \\
  &= \underset{s=0} \to {\text{res}}\,\,  \check g(s)R_K(s+1) + \frac{1}{2\pi}
\int_{(-\frac{1}{2})} \check g(s) R_K(s+1) ds.
\endalign
$$
Hence using the above estimates we get
$$
G = \int^\infty_0 g(y)(\alpha \log y + \beta) \frac{dy}{y} +
O((\frac{q}{x})^{\frac{1}{2}}). \tag6.48
$$
\enddemo
\bigskip

\proclaim{Corollary 6.3}  For $Y \geqslant 2X \geqslant 2$ we have
$$
\dsize \sum_{X \leqslant n \leqslant Y} \tau^2(n,\chi)n^{-1} \ll {\Cal L}(Y)
\log\frac{Y}{X} + (\frac{q}{X})^{\frac{1}{2}} \tag6.49
$$
where
$$
{\Cal L}(Y) = L(1,\chi) \bigl( L(1,\chi)\log Y + |L'(1,\chi)| \bigr). \tag6.50
$$
\endproclaim
\medskip
For $G$ given by (6.40) the formula (6.48) becomes
$$
G = \int^\infty_0 |a(y)|^2 (\alpha \log y + \beta) \frac{dy}{y} +
O((\frac{q}{T})^{\frac{1}{2}}). \tag6.51
$$
Hence, if $2T \leqslant Y \leqslant T^{\frac{31}{30}}$ we get
$$
G \ll {\Cal L}(T) \log \frac{Y}{T} + (\frac{q}{T})^{\frac{1}{2}}.
$$
\mn
Introducing this into (6.35) we end up with the following
\proclaim{Proposition 6.4}  Let $\lambda(n)$ be the coefficients of an
automorphic form given by Hecke characters of the imaginary quadratic
field $K = \Bbb Q(\sqrt{-q})$ of discriminant $-q$.  Let $a(y)$ be a
function satisfying (6.38) with $Y = qT$ and $T \geqslant K^{32}
q^{65}$. Then
$$
\int^{2T}_T |\dsize \sum_n a(n) \lambda(n)n^{-\frac{1}{2} -it}|^2 dt
\ll T {\Cal L}(T) \log q \tag6.52
$$
where ${\Cal L}(T)$ is defined by (6.50) and the implied constant is absolute.
\endproclaim
\bigskip

\remark{Remark}  The above bound comes from the diagonal terms.  The
other terms contribute slightly less, namely $T {\Cal L}(q)\log q$.
\endremark
\bn
\subhead{7. Approximate Functional Equation} \endsubhead
\medskip
We restrict our attention to $L$-functions for class group
characters of an imaginary quadratic field $K = \Bbb Q(\sqrt{-q})$
where $-q$ is the discriminant.  We assume that $q$ is odd and $q
> 4$, so $q \equiv 3(\text{mod } 4)$ and $q$ is squarefree. Fix
$\psi \in {\widehat{\Cal C}\ell}(K)$ and put $$ \lambda(n) =
\dsize \sum_{N{\frak a} = n} \psi({\frak a}). \tag7.1 $$ These are
Hecke eigenvalues of an automorphic form (a theta series) of
weight $k=1$, level $q$ and character  $$ \chi(n) = \bigl(
\frac{n}{q} \bigr), \tag7.2  $$ the Jacobi symbol. Let $$ L(s) =
\dsize \sum^\infty_{n=1} \lambda(n)n^{-s} \tag7.3 $$ be the
corresponding Hecke $L$-function.  For example, if $\psi$ is a
genus character then $$ L(s) = L(s,\chi_v)L(s,\chi_w) \tag7.4 $$
where $\chi_v,\chi_w$ are the real characters of conductor $v,w$
respectively with $vw = q$, i.e. $$ \chi_v(n) = \bigl( \frac{n}{v}
\bigr),  \qquad \chi_w(n) = \bigl( \frac{n}{w} \bigr) \tag7.5 $$
are the corresponding Jacobi symbols.  Observe that
$\chi_v,\chi_w$ are characters for real and imaginary quadratic
fields. If $\psi \in {\widehat{\Cal C}\ell}(K)$ is not a genus
character (i.e. $\psi$ is not real) then the corresponding
$L$-function does not factor into Dirichlet $L$-function. However,
in any case the complete product $$ \Lambda(s) = Q^s \Gamma(s)L(s)
\quad \text{ with } \quad Q = \frac{\sqrt q}{2 \pi} \tag7.6 $$ has
analytic continuation to the whole complex $s$-plane, except for a
simple pole at $s=1$ if $\psi$ is the trivial character, in which
case $$ L(s) = \zeta_K(s) = \zeta(s)L(s,\chi) \tag7.7 $$ is the
zeta function of $K$. Moreover for any $\psi \in {\widehat{\Cal
C}\ell}(K)$ we have the functional equation (which is due to Hecke
, see also (3.9)) $$ \Lambda(s) = \Lambda(1-s). \tag7.8 $$

In this section we derive a Dirichlet series representation of
$L(s)$ tempered by a test function which makes the series rapidly
convergent. Formulas of this type are known in the literature as
``approximate functional equations".  In our context, this is a
somewhat misleading name, because we need exact expressions to be
able to differentiate.  We rather think of these as a kind of
Poisson's summation formulas.

Let $G(u)$ be a holomorphic function in the strip $|\text{Re }u|
\leqslant 1$ such that
$$
\align
G(u) = G(-u) \tag7.9 \\
G(0) = 1 \tag7.10 \\
G(u) \ll 1. \tag7.11
\endalign
$$ \mn Consider the integral $$ I(s) = \frac{1}{2 \pi i}
\int_{(1)} \Lambda(s+u)G(u) u^{-1} du $$ for $0 < \text{ Re } s <
1$. Moving the path of integration to the line Re $u = -1$ and
applying (7.8) we get $$ \Lambda(s) = I(s) + I(1-s) -
\frac{G(s-1)}{s-1} \underset{u=1}\to {\text{res}} \Lambda(u). $$
On the other hand, introducing the Dirichlet series (7.3) and
integrating termwise we obtain $$ I(s) = \dsize \sum_n \lambda(n)
\frac{1}{2 \pi i} \int_{(1)} \bigl( \frac{Q}{n}\bigr)^{s+u}
\Gamma(s+u)G(u)u^{-1} du. $$ From both expressions we obtain
(after dividing by $(\frac{\sqrt q}{2 \pi})^s \Gamma(s))$.
\bigskip

\proclaim{Proposition 7.1}  For $s$ with $0 < \text{ Re } s < 1$ we have
$$
\align
L(s) = \dsize \sum_n \lambda(n) n^{-s} V_s(\frac{n}{Q}) &+ X(s) \dsize
\sum_n \lambda(n)n^{s-1} V_{1-s}(\frac{n}{Q}) \tag 7.12 \\
  &- \frac{G(s-1)}{(s-1)\Gamma(s)} Q^{1-s} L(1,\chi),
\endalign
$$
where
$$
X(s) = Q^{1-2s} \Gamma(1-s)/\Gamma(s) \tag7.13
$$
$$
V_s(y) = \frac{1}{2 \pi i} \int_{(1)} \frac{\Gamma(s+u)}{\Gamma(s)} \,
\frac{G(u)}{u} y^{-u} du \tag7.14
$$
and the last (the residual) term in (7.12) exists only if $\psi$ is
the trivial character of ${\Cal C}\ell(K)$.
\endproclaim
\medskip
We shall apply (7.12) for points on the critical line Re $s =
\frac{1}{2}$.  Choosing
$$
G(u) = (\text{cos } \frac{\pi u}{A})^{-A} \tag7.15
$$
where $A \geqslant 4$ is a fixed integer we derive
\medskip

\proclaim{Lemma 7.2}  If Re $s = \frac{1}{2}$ then
$$
y^aV^{(a)}_s(y) = \delta(a) - \frac{G(s)}{\Gamma(s+1-a)} y^s +
O(\frac{y}{|s|}) \tag7.16
$$
$$
y^a V^{(a)}_s(y) \ll (1 + \frac{y}{|s|})^{-A} \tag7.17
$$
\mn
for any $a \geqslant 0$, the implied constant depending only an $a$ and $A$
(here $\delta(0) =1$ and $\delta(a) =0$ if $a > 0$).
\endproclaim
\bigskip

\remark{Remark}  We have $G(s) \ll e^{- \pi|s|}$ and
$\Gamma(s+1-a)^{-1} \ll |s|^{a-1} e^{\frac{\pi}{2}|s|}$ ;hence,
Lemma 7.2 yields $$ y^a V^{(a)}_s(y) = \delta(a) + O \biggl(
\sqrt{\frac{y}{|s|}} \biggr). $$
\endremark
\bigskip

\demo{Proof}  Differentiating (7.14) $a$ times we get
$$
y^a V^{(a)}_s(y) = \frac{1}{2 \pi i} \int_{(1)}
\frac{\Gamma(s+u)}{\Gamma(s)} \, \frac{G(u)}{u} \,
\dsize \prod_{0 \leqslant b < a} (-u-b)\,\,\, y^{-u} du. \tag7.18
$$
For the proof of (7.16) we move the integration to the line Re $u =
-1$ getting the first two terms as residues at $u=0$ and $u = -s$
respectively.  Using Stirling's formula the resulting integral
on Re $u = -1$ is estimated by
$$
\align
\int_{(-1)} &|s+u|^{-1} e^{-\frac{\pi}{2}|s+u|+\frac{\pi}{2}|s|-\pi|u|}
|u|^{a-1} y|du| \\
  &\ll \int_{(-1)} |s+u|^{-1} e^{-\frac{\pi}{2}|u|} |u|^{a-1} y|du|
\ll \frac{y}{|s|}.
\endalign
$$
For the proof of (7.17) we move the integration to the line Re $u =
A$.  Using Stirling's formula the resulting integral is estimated by
$$
\int_{(A)} |s+u|^A e^{-\frac{\pi}{2}|u|} |u|^{a-1} y^{-A} |du|
\ll \bigl( \frac{|s|}{y} \bigr)^A.
$$
This yields (7.17) if $y > |s|$.  In the case $y \leqslant |s|$ we get
(7.17) from (7.16).

Actually we shall apply (7.12) to estimate the quotients $$
\ell(s) = \frac{L(s)-L(s')}{s-s'} \tag7.19 $$ for points $s,s'$ on
the critical line (if $s=s'$, then $\ell(s) = L'(s)$ is the
derivative of $L(s)$).  Here we do not display the dependence of
$\ell(s)$ on the second point $s'$ for notational simplicity. This
abbreviated notation (also used for other forthcoming quotients)
will be justified when we fix $s'$ in terms of $s$. Put $$ \align
x(s) = \frac{X(s)-X(s')}{s-s'}, \tag7.20 \\ v_s(y) =
\frac{V_s(y)-V_{s'}(y)}{s-s'}, \tag7.21 \\ w_s(y) =
\frac{1-y^{s-s'}}{s-s'}. \tag7.22
\endalign
$$
From (7.12) we derive (by adding and subtracting terms)
$$
\align
\ell(s) &= \biggl( \dsize \sum_n - X(s) \overline{\dsize \sum_n} \biggr)
\lambda(n)n^{-s}
w_s(n) V_s (\frac{n}{Q}) \tag7.23 \\
  &+ \biggl( \dsize \sum_n - X(s) \overline{\dsize \sum_n} \biggr) \lambda(n)n^{-s}
 v_s (\frac{n}{Q}) \\
  &+ x(s) \overline{\dsize \sum_n} \lambda(n)n^{-s} V_s(\frac{n}{Q}) + O(\frac{n}{|s|}),
\endalign
$$
where $\underset n\to {\bar \Sigma}$ stands for the complex conjugate of
$\underset n\to \Sigma$.

Now we need estimates for $x(s)$ and for derivatives of $v_s(y),w_s(y)$.
\enddemo
\bigskip

\proclaim{Lemma 7.3}  For $s,s'$ on the critical line we have $$
|w_s(y)| \leqslant |\log y|,  \qquad |w'_s(y)| = y^{-1}. \tag7.24
$$
\endproclaim
\bigskip

\proclaim{Lemma 7.4}  For $s,s'$ on the critical line we have
$$
y^a v^{(a)}_s(y) \ll \bigl(\frac{y}{|s|}\bigr)^{\frac{1}{4}} \bigl(1 +
\frac{y}{|s|}\bigr)^{-A} \tag7.25
$$
if $a \geqslant 0$, the implied constant depends only on $a$ and $A$.
\endproclaim
\bigskip

\demo{Proof}  If $|s-s'| > 1$ then Lemma 7.4 follows from Lemma 7.2 by
subtracting the estimates.  Let $|s-s'| \leqslant 1$.  Subtract (7.18) for
$s'$ from that for $s$ and divide by $s-s'$ to obtain a corresponding
expression for derivatives of $v_s(y)$.  Then move the integration
from Re $u=1$ to Re $u= \alpha$
with $-\frac{1}{4} \leqslant \alpha \leqslant A$.
Note that there is no pole at $u=0$.  Then estimate as follows
$$
\align
\hskip40pt \frac{1}{|s-s'|} \biggl| \frac{\Gamma(s+u)}{\Gamma(s)} &-
\frac{\Gamma(s'+u)}{\Gamma(s')} \biggr| \leqslant \biggl|
\frac{\partial}{\partial s} \, \frac{\Gamma(s+u)}{\Gamma(s)}
\biggr|_{s=s_0} \\
 &= \biggl| \frac{\Gamma(s_0+u)}{\Gamma(s_0)} \biggr|
|\psi(s_0+u)-\psi(s_0)| \ll |s_0 +u|^\alpha \,
e^{\frac{\pi}{2}|u|}|u|
\endalign
$$
where $s_0$ is a point on the critical line between $s$ and $s'$.
Moreover $G(u) \ll e^{-\pi|u|}$.  Hence
$$
y^a v^{(a)}_s(y) \ll y^{-\alpha} \int_{(\alpha)} |s_0+u|^\alpha |u|^a
\, e^{-\frac{\pi}{2}|u|} |du| \ll \bigl(\frac{|s|}{y}\bigr)^\alpha.
$$
This implies (7.25) by taking $\alpha = -\frac{1}{4}$ if $y \leqslant |s|$,
or $\alpha = A$ if $y > |s|$.

Before estimating $x(s)$ note that $|X(s)|=1$ for Re $s =
\frac{1}{2}$ ; more precisely $$ X(\tfrac{1}{2} + it) = Q^{-2\,
it}\, \frac{\Gamma(\frac{1}{2}-it)}{\Gamma(\frac{1}{2} +it)} =
(\frac{e}{tQ})^{2it}\, \{1 + \varepsilon(t)\} \tag7.26 $$ if $t
\geqslant 1$, where $\varepsilon(t) \ll t^{-1}$ and
$\varepsilon'(t) \ll t^{-2}$.  Hence we derive
\enddemo
\bigskip

\proclaim{Lemma 7.5}  For $s = \frac{1}{2} + it$ and $s =
\frac{1}{2} + it'$ with $t,t' \geqslant 1$ we have $$ x(s) = -2\,
(\frac{e}{tQ})^{it}\,(\frac{e}{t'Q})^{it'}\, \, \frac{
\sin(t-t')\log tQ}{t-t'} + O(\frac{1}{t}); \tag7.27 $$
consequently, $$ |x(s)| = 2\, \bigl| \frac{\sin(t-t')\log
tQ}{t-t'} \bigr| + O(\frac{1}{t}). \tag7.28 $$
\endproclaim
\bigskip

\demo{Proof}  If $|t-t'| < \frac{t}{2}$ then (7.27) follows from
(7.26); otherwise (7.27) is trivial.
\enddemo
\bigskip

Applying the inequality sin $x \geqslant \alpha x$ if $0 \leqslant x \leqslant
\pi(1-\alpha)$ we get
\proclaim{Corollary 7.6}  Let
$0 \leqslant \alpha \leqslant 1$.  If \,$|t-t'|\log tQ
\leqslant \pi(1-\alpha)$, then
$$
|x(s)| \geqslant 2 \alpha \log \,tQ + O(\frac{1}{t}). \tag7.29
$$
\endproclaim
\mn
\subhead{8. Evaluation of $\ell(s)$ on Average} \endsubhead
\medskip
Our goal is to eliminate most of the terms in (7.23) by estimating
them on average with respect to a well-spaced set of points $s$ on
the critical line.  We begin by any set, say $S(T)$, of points $$
s_r = \tfrac{1}{2} + it_r,\,\,\,\, r = 1,2,\dotsc,R \tag8.1 $$
such that for $T \geqslant 2$ $$ T < t_1 < t_2 < \ldots < t_R
\leqslant 2T, \tag8.2 $$  $$ t_{r+1} - t_r \geqslant 1,\,\, \text{
if } 1 \leqslant r < R. \tag8.3 $$ To each point $s_r$ we
associate a point $$ s'_r = \tfrac{1}{2} + it'_r. \tag8.4 $$
\bigskip

\remark{Remarks}  The companion $s'_r$ to $s_r$ may not be in
$S(T)$. Actually our main interest will be to choose $s'_r$ very
close to $s_r$.  For example $s_r,s'_r$ can be consecutive zeros
of $L(s)$ on the critical line.  We may have $s_r = s'_r$ if this
is a double zero. For the time being we assume that $T \geqslant
q^{65}$ to comply with the condition of Proposition 6.4, but after
shaping the  basic estimates this assumption can be dispensed
because the results hold true trivially if $T < q^{65}$.

First we estimate the sum $$ A_1(s) = \dsize \sum_n \lambda(n)
n^{-s} w_s(n)V_s(\frac{n}{Q}) \tag8.5 $$ \mn on average with
respect to the points $s \in S(T)$.  Recall that $w_s(y)$
satisfies (7.24) and $V_s(y)$ satisfies (7.17).  We partition this
sum smoothly into three sums, say $A_1(s) = A_{11}(s) + A_{12}(s)
+ A_{13}(s)$, where the partial sums are supported on the segments
$n_1 \ll q^4 \ll n_2 \ll T \ll n_3$, respectively.  For estimation
of $A_{11}(s)$ we apply Lemma 5.3 with $c \ll \log q$ and $a_n \ll
|\lambda(n)|n^{-\frac{1}{2}} \log n$ and obtain $$ \dsize \sum_s
|A_{11}(s)|^2 \ll T(\dsize \sum_{n \ll q^4} \tau^2(n)n^{-1} \log^2
n)\log q \ll T (\log q)^7 $$ by the trivial estimate (6.15).  For
estimation of $A_{12}(s)$ we apply Lemma 5.3 with $c \ll \log T$
and $a_n \ll |\lambda(n)|n^{-\frac{1}{2}} \log n$ ; now, however,
we take advantage of the better bound given by  (6.42) to see that
$$ \dsize \sum_s |A_{12}(s)|^2 \ll T(\dsize \sum_{q^4 \ll n \ll T}
\tau^2(n,\chi)n^{-1} \log^2 n)\log T \ll T {\Cal L}(T)(\log T)^4
$$ by (6.50).  For estimation of $A_{13}(s)$ we apply Proposition
5.4 with $c \ll \log q$ and $a_n \ll |\lambda(n)|n^{-\frac{1}{2}}
\log n$, together with Proposition 6.4 getting $$ \dsize \sum_s
|A_{13}(s)|^2 \ll T {\Cal L}(T)(\log T)^2(\log q)^2. $$ Next we
estimate the sum $$ A_2(s) = \dsize \sum_n \lambda(n)n^{-s}
v_s(\frac{n}{Q}). \tag8.6 $$ \mn The arguments are the same as
those applied for $A_1(s)$ above, and the corresponding estimates
are sharper by two logarithms because $v_s(y) \ll 1$ and $w_s(y)
\ll \log y$.  Precisely, we get $A_2(s) = A_{21}(s) + A_{22}(s) +
A_{23}(s)$ with $$ \dsize \sum_s |A_{21}(s)|^2 \ll T(\log q)^5 $$
$$ \dsize \sum_s |A_{22}(s)|^2 \ll T{\Cal L}(T)(\log T)^2 $$ $$
\dsize \sum_s |A_{23}(s)|^2 \ll T{\Cal L}(T)(\log q)^2. $$
\endremark
\bigskip

\remark{Remark}  If we used the more precise bound for $v_s(y)$
given in Lemma 7.4, then the above estimates could be improved
further, but that leads to no advantage here.

It remains to estimate the sum $$ A_3(s) = \dsize \sum_n
\lambda(n)n^{-s} V_s(\frac{n}{Q}). \tag8.7 $$ Similarly we
partition this sum smoothly into $A_3(s) = A_{31}(s) + A_{32}(s) +
A_{33}(s)$, and apply the same arguments as those for $A_1(s)$,
getting $$ \dsize \sum_s |A_{31}(s)|^2 \ll T(\log q)^5 $$ $$
\dsize \sum_s |A_{32}(s)|^2 \ll T{\Cal L}(T)(\log T)^2 $$ $$
\dsize \sum_s |A_{33}(s)|^2 \ll T{\Cal L}(T)(\log q)^2. $$
However, we are not satisfied with the above bound for
$A_{31}(s)$. First we clear from $A_{31}(s)$ the factor $$
V_s(\frac{n}{Q}) = 1 + O \biggl( \sqrt{\frac{n}{QT}} \biggr) $$
(see Lemma 7.2) and replace the smooth cut-off function (from the
relevant partition) in the range $n \asymp q^4$.  We get $$
A_{31}(s) = N(s) + \tilde N(s) + O(q^4 T^{-\frac{1}{2}}) $$ where
$$ N(s) = \dsize \sum_{n \leqslant q^4} \lambda(n) n^{-s} \tag8.8
$$  $$ \tilde N(s) = \dsize \sum_{x < n \leqslant y} \lambda(n)
\alpha(n) n^{-s} \tag8.9 $$ for some $x,y$ with $q^4 \ll x < y \ll
q^4$ and $\alpha(n) \ll 1$.  By Lemma 5.3 and Corollary 6.3 we
derive $$ \dsize \sum_s |\tilde N(s)|^2 \ll T \dsize \sum_{x < x
\leqslant y} \tau^2(n,\chi)n^{-1} \ll T{\Cal L}(q)\log q. $$
Moreover the error term $O(q^4 T^{-\frac{1}{2}})$ contributes at
most $R(q^4 T^{-\frac{1}{2}})^2 \ll q^8$, which is absorbed by
$T{\Cal L}(q) \log q$.  Therefore we have $$ \dsize \sum_s
|A_{31}(s) - N(s)|^2 \ll T{\Cal L}(q)\log q. $$
\endremark
\bigskip

Gathering the above estimates together with (7.23) we obtain
\proclaim{Proposition 8.1}  Let $S(T)$ be a set of points satisfying
(8.1)-(8.3) with $T \geqslant 2$.
Put
$$
D(T) = \dsize \sum_s |\ell(s)-x(s)\overline N(s)|^2
$$
where $s$ runs over $S(T)$
(recall the settings (7.19), (7.20), (8.8)).  We have
$$
D(T) \ll T(\log q)^7 + T{\Cal L}(T)(\log T)^4 \tag8.10
$$
where ${\Cal L}(T)$ is defined by (6.50), the implied constant being
absolute.
\endproclaim
\medskip
Assuming that $L(1,\chi)$ is small relatively to log $T$ (so is ${\Cal
L}(T)$) we can interpret the bound (8.10) as saying that $x(s)\bar
N(s)$ approximates to $\ell(s)$ at almost all points $s$ in any
well-spaced set $S(T)$.
\bn
\subhead{9. Estimation of $x(s)$ on Average} \endsubhead
\medskip
Recall that the Hecke $L$-function for a character $\psi \in
{\widehat{\Cal C}\ell}(K)$ has the Euler product $$ L(s) = \dsize
\prod_{\frak p} (1-\psi({\frak p})(N{\frak p})^{-s})^{-1} = \dsize
\sum_n \lambda(n)n^{-s}; \tag9.1 $$ similarly the inverse
satisfies $$ L^{-1}(s) = \dsize \prod_{\frak p} (1-\psi({\frak
p})(N{\frak p})^{-s}) = \dsize \sum_m \lambda^*(m)m^{-s}, \tag9.2
$$ say, where $$ \lambda^*(m) = \dsize \sum_{N{\frak a} =m}
\mu({\frak a}) \psi({\frak a}).\tag9.3 $$ Note that
$\lambda^*(m)$, like $\lambda(m)$, often vanishes if the class
number is small. We have

$$ |\lambda^*(m)| \leqslant \tau(m,\chi). \tag9.4 $$ Hence we have
a reason to believe that the partial sum of $L^{-1}(s)$ $$ M(s) =
\dsize \sum_{m \leqslant q^4} \lambda^*(m)m^{-s} \tag9.5 $$
approximates to $N^{-1}(s)$ at almost all points $s$ on the
critical line. Our goal is to estimate the sum $$ E(T) = \dsize
\sum_s |\ell(s) {\overline M}(s) - x(s)|. \tag9.6 $$

We begin by writing $M(s)N(s) = 1 + B(s)$, where $$ B(s) = \dsize
\sum_{q^4 < \ell \leqslant q^8} b(\ell)\ell^{-s} $$ with $$
b(\ell) = \dsize \sum \Sb mn=\ell \\ m,n \leqslant q^4 \endSb
\lambda^*(m)\lambda(n). $$ Then we split $E(T)$ as follows : $$
\align E(T) &= \dsize \sum_s |(\ell(s)-x(s){\overline
N}(s)){\overline M}(s) + x(s){\overline B}(s)| \\
  &\ll D(T)^{\frac{1}{2}}(\dsize \sum_s |M(s)|^2)^{\frac{1}{2}} +
(\log T)\dsize \sum_s |B(s)|.
\endalign
$$ Here we have $$ \dsize \sum_s |M(s)|^2 \ll T(\log q) \dsize
\sum_{m \leqslant q^4} \tau^2(m)m^{-1} \ll T(\log q)^5 $$ and $$
\dsize \sum_s |B(s)| = \dsize \sum_s |\dsize \sum \Sb m,n\leqslant
q^4 \\ mn > q^4 \endSb \lambda^*(m)\lambda(n)(mn)^{-s}|. $$ Note
the condition $mn > q^4$ implies that either $m$ or $n$ is larger
than $q^2$.  Having this information recorded we relax the
condition $mn > q^4$ by any method of separation of variables, for
example by applying Lemma 9 of [DFI1]. This separation costs us a
factor $\log q$. It follows that $$ \dsize \sum_s |B(s)| \ll
T(\log q)^2 (\dsize \sum_{q^2 < m \leqslant q^4}
\tau^2(m,\chi)m^{-1})^{\frac{1}{2}}(\dsize \sum_{n \leqslant q^4}
\tau^2(n)n^{-1})^{\frac{1}{2}}. $$ By (6.49) we derive $$ \dsize
\sum_s |B(s)| \ll T{\Cal L}(q)^{\frac{1}{2}} (\log
q)^{\frac{9}{2}}. $$ \bn These estimates yield
\proclaim{Proposition 9.1}  Let $S(T)$ be a set of points
satisfying (8.1) - (8.3) with $T \geqslant 2$.  Then $$ E(T) \ll
T(\log q)^6 + T{\Cal L}(T)^{\frac{1}{2}}(\log T)^2(\log
q)^{\frac{5}{2}} \tag9.7 $$ where the implied constant is
absolute.
\endproclaim
\medskip
Assuming that $L(1,\chi)$ is relatively small Proposition 9.1
asserts that $\ell(s)\bar M(s)$ approximates to $x(s)$ at almost
all points $s$ in any well-spaced set $S(T)$.  This assertion is
particularly interesting if $\ell(s) = (L(s)-L(s'))(s-s')^{-1}$ is
very small, because it implies that $x(s) =
(X(s)-X(s'))(s-s')^{-1}$ is also quite small. Put $$ \Delta(T) =
\dsize \sum_s |\ell(s)|^2 = \dsize \sum_s \biggl|
\frac{L(s)-L(s')}{s-s'} \biggr|^2. \tag9.8 $$ By Cauchy's
inequality we get $$ \dsize \sum_s |\ell(s)M(s)| \ll (T
\Delta(T)\log^5 q)^{\frac{1}{2}}. $$ Applying this to $E(T)$ in
(9.6) we derive by (9.7) the following estimate $$ \dsize \sum_s
|x(s)| \ll T(\log q)^6 + T{\Cal L}(T)^{\frac{1}{2}}(\log T)^2(\log
q)^{\frac{5}{2}} + (T \Delta(T)\log^5 q)^{\frac{1}{2}}. \tag9.9 $$
We shall make the estimate (9.9) more explicit by cosmetic
preparations. First on the left-hand side we use (see (7.28)) $$
|x(s)| = 2\, \biggl| \frac{\text{sin}(t-t')\log t}{t-t'} \biggr| +
O(\log q). \tag9.10 $$ Note that the error term $O(\log q)$
contributes in total $O(R \log q)$ which is absorbed by the first
term $T(\log q)^6$ on the right-hand side of (9.9).  Next we
replace ${\Cal L}(T)$ in (9.9) by $L(1,\chi)\log q$.  This can be
justified, because the modified inequality is trivial unless $$
(\log T)L(1,\chi)^{\frac{1}{2}}(\log q)^3 \leqslant 1. \tag9.11 $$
Moreover, if (9.11) holds then we find that ${\Cal L}(T) \ll
L(1,\chi)\log q$.  Finally, we no longer restrict the points $s =
\frac{1}{2}+it$ to a dyadic segment $T < t \leqslant 2T$.  The
extension to the segment $1 \leqslant t \leqslant T$ can be now
derived by adding the new inequalities (9.9) for sets of points in
the segments $2^\nu \leqslant t \leqslant 2^{\nu+1}$ with $1
\leqslant 2^\nu \leqslant T$.  We state the result in a
self-contained format.
\bigskip

\proclaim{Proposition 9.2}  Let $s$ run over a set of points on the
critical line $s = \frac{1}{2} +it$ with $2 \leqslant t \leqslant T$ which are
spaced by at least one.  To every $s$ in the set we associate a point
$s' = \frac{1}{2} + it'$.  Then we have
$$
\align
\dsize \sum_s \biggl| \frac{\text{\rm{sin}\/}(t-t')\log t}
{(t-t')\log t} \biggr| &\ll
\frac{T}{\log T} (\log q)^6 + T(\log T)L(1,\chi)^{\frac{1}{2}}(\log q)^3
\tag9.12 \\
  &+\frac{(\log q)^{\frac{5}{2}}}{\log T} \biggl( T \dsize \sum_s \biggl| \frac{L(s)-L(s')}{s-s'} \biggr|^2
\biggr)^{\frac{1}{2}} .
\endalign
$$
\endproclaim
\medskip
This is our principal estimate from which one can deduce numerous
attractive propositions.  But first we wish to emphasize that
(9.12) has no permanent value; it has some quality only in the
absence of the Riemann hypothesis.  Indeed, assuming only the
lower bound $$ L(1,\chi) \gg (\log q)^{-6} \tag9.13 $$ (recall
that the Riemann hypothesis for $L(s,\chi)$ yields (1.5)) we find
that the middle term on the right side of (9.12) is bounded below
by $T{\log T}$.  On the other hand the left side of (9.12) is
trivially bounded by $R \leqslant T$. Therefore our principal
estimate (9.12) is insignificant if (9.13) is true.  We certainly
believe in the truth of (9.13), nevertheless as long as
$L(1,\chi)$ is not proved to be relatively large (the best known
unconditional estimate being $L(1,\chi) \gg
q^{-\varepsilon}$,which is not effective), there are some valuable
features of (9.12). \bn \subhead{10. Applications} \endsubhead
\medskip
In this section we derive a few consequences of the principal estimate
(9.12). We begin by eliminating the last term $(T
\Delta(T))^{\frac{1}{2}}(\log q)^{\frac{5}{2}}(\log T)^{-1}$.

    If all the points $s$ and their companions $s'$ are zeros of $L(s)$
(double zeros if $s = s'$) on the critical line, then $$ \ell(s) =
\frac{L(s)-L(s')}{s-s'} = 0, \tag10.1 $$ and consequently
$\Delta(T)=0$.  Actually we do not require $s$ and $s'$ to be
zeros of $L(s)$; the condition (10.1) means that $s$ and its
companion $s'$ are on the same level curve of $L(s)$ (and
$L'(s)=0$ if $s=s'$).  We can still assume less than (10.1). For
example if $s$ and its companion $s'$ satisfy $$ \bigl|
\frac{L(s)-L(s')}{s-s'} \bigr| \leqslant (\log q)^{\frac{7}{2}}
\tag10.2 $$ then $\Delta(T) \leqslant T(\log q)^7$, so on the
right side of (9.12) the last term is absorbed by the first one.
From now on we assume that the points $s$ and their companions
$s'$ satisfy (10.2).  For the points so chosen the estimate (9.12)
reduces to $$ \dsize \sum_s \bigl| \frac{\text{sin}(t-t')\log
t}{(t-t')\log t} \bigr| \ll \frac{T}{\log T} (\log q)^6 + T(\log
T)L(1,\chi)^{\frac{1}{2}}(\log q)^3. \tag10.3 $$ Choose any $T$
with $$ (\log q)^{A+6} \leqslant \log T \leqslant
L(1,\chi)^{-\frac{1}{2}}(\log q)^{-A-3}. \tag10.4 $$ where $A$ is
a positive constant; then  (10.3) implies that $$ \dsize \sum_s
\bigl| \frac{\text{sin}(t-t')\log t}{(t-t')\log t} \bigr| \ll
T(\log q)^{-A}. \tag10.5 $$ We tacitly assumed that $L(1,\chi)$ is
small to be sure that the interval (10.4) is not void ; precisely
for this reason we require $$ L(1,\chi) \leqslant (\log
q)^{-4A-18}. \tag10.6 $$ Let $R = R(\alpha,T)$ be the number of
points $s$ in the set with the companions $s'$ satisfying (10.2)
and $$ |t-t'| \leqslant \frac{\pi(1-\alpha)}{\log t} \tag10.7 $$
where $0 < \alpha \leqslant 1$.  For such points we have $$
\frac{\text{sin}(t-t')\log t}{(t-t')\log t} \geqslant \alpha.
\tag10.8 $$ Hence (10.5) gives the following bound for the number
of points in question : $$ R \ll \alpha^{-1} T(\log q)^{-A}
\tag10.9 $$ where the implied constant is absolute.

In conclusion we rephrase the obtained results in a positive mode.
\bigskip

\proclaim{Proposition 10.1}  Let $A \geqslant 0$ and $\log T
\geqslant (\log q)^{A+6}$. Suppose there is a set of points $$
S(T) = \{s_r = \tfrac{1}{2} + it_r;\,\,\,2 \leqslant t_1 < \ldots
< t_R \leqslant T,\,\,\,t_{r+1} - t_r \geqslant 1\} $$ and a set
of companions $S'(T) = \{s'_r = \frac{1}{2} + it'_r ;\,\,r =
1,\dotsc,R\}$ such that $$  |t_r-t'_r| \leqslant \frac{\pi(1-
\alpha)}{\log {t_r}},\,\, \text{ with } 0 < \alpha \leqslant 1,
\tag10.10 $$ $$ \bigl| \frac{L(s_r)-L(s'_r)}{s_r-s'_r} \bigr|
\leqslant (\log q)^{\frac{7}{2}}. \tag10.11 $$
 Suppose the number of points in the set satisfies $$ R = |S(T)|
\geqslant \frac{cT}{\alpha(\log q)^A} \tag10.12 $$ where $c$ is a
large absolute constant, effectively computable.  Then $$
L(1,\chi) \geqslant (\log T)^{-2}(\log q)^{-2A-6}. \tag10.13 $$
\endproclaim
\medskip
We certainly believe that any Hecke $L$-function satisfies the
conditions of Proposition 10.1, provided $q$ is large.
We recommend the points $s_r$ to be
zeros of $L(s)$ and $s'_r$ to be the nearest zero to $s_r$ on the
critical line (if $s_r$ has order two or more, then $s'_r =s_r$).  For
this choice (10.11) holds automatically, while (10.10) asserts that
the gaps between chosen pairs of zeros is smaller than the normal
average spacing.  We need a considerable number of such small gaps
between consecutive zeros, but less than the true order of magnitude.
In particular taking $A=6$ and $\log T = (\log q)^{12}$ we get
$L(1,\chi) \geqslant (\log q)^{-42}$,
provided the number of well-spaced zeros of sub-normal gaps and height
up to $T$ is at least $T(\log T)^{-\frac{1}{2}}$.
\medskip
\remark{Remark}  If the points of $S(T)$ are zeros of $L(s)$ then
the condition that they are spaced by at least one can be dropped
at the expense of extra factor log $T$ in (10.12) (see how this is
justified in the proof of Corollary 10.2).  Hence we derive
Theorem 1.1.
\endremark
\medskip

An interesting case is the trivial class group character $\psi=1$.  In
this case the Riemann zeta function appears as a factor of the Hecke
$L$-function, $L(s) = \zeta_K(s) = \zeta(s)L(s,\chi)$, so we can
choose the zeros of $L(s)$ from those of $\zeta(s)$ and state the
conditions without ever mentioning the exceptional conductor $q$.
Note we have precise control over $\alpha$.
Taking $A=12$ and $\alpha = (\log T)^{-\frac{1}{2}}$ we derive
from Proposition 10.1 the following
\medskip

\proclaim{Corollary 10.2}  Let $\rho = \frac{1}{2} + i\gamma$
denote the zeros of $\zeta(s)$ on the critical line and $\rho' =
\frac{1}{2} +i\gamma'$ denote the nearest zero to $\rho$ on the
critical line ($\rho' = \rho$ if it is a multiple zero).  Suppose
that $$ \# \biggl\{ \rho;\,\,\,0 < j \leqslant
T,\,\,\,|\gamma-\gamma'| \leqslant \frac{\pi}{\log \gamma} \,
\bigl( 1 - \frac{1}{\sqrt{\log \gamma}} \bigr) \biggr\} \gg T(\log
T)^{\frac{4}{5}} \tag10.14 $$ for any $T \geqslant 2001$.  Then we
have $$ L(1,\chi) \gg (\log q)^{-90} \tag10.15 $$ where the
implied constant is effectively computable in terms of that in
(10.14).
\endproclaim
\bigskip

\demo{Proof}  The number of zeros $\rho = \frac{1}{2} +i\gamma$ with
$t < \gamma \leqslant t+1$ is bounded
by $O(\log t)$.  Therefore one can select from the
set of zeros in (10.14) a subset of well-spaced points of cardinality
$R \gg T(\log T)^{-\frac{1}{5}}$.  This subset satisfies the conditions
of Proposition 10.1 with $A=14,\,\log T = (\log q)^{20}$ and $\alpha =
(\log T)^{-\frac{1}{2}} = (\log q)^{-10}$, provided $q$ is sufficiently
large, giving $L(1,\chi) \geqslant (\log q)^{-90}$.  For small $q$ the lower
bound (10.15) is obtained by adjusting the implied constant.
\enddemo
\bigskip

\remark{Remark}  The normal average gap between consecutive zeros
$\rho = \frac{1}{2} +i \gamma,\rho' = \frac{1}{2} +i \gamma'$ of $\zeta(s)$ is
$2\pi(\log \gamma)^{-1}$.  Hence the condition (10.14) refers to gaps
slightly smaller than the half of the average.
\endremark

\enddocument
\bye